\documentclass{amsart}

\usepackage[bb=dsserif]{mathalpha}
\usepackage{bm}

\usepackage{graphicx} 
\usepackage{geometry}
\usepackage[foot]{amsaddr}
\usepackage[english]{babel}
\usepackage{amsthm}

\usepackage{url}
\usepackage{breakurl}
\usepackage[breaklinks]{hyperref}
\usepackage{xfrac}
 
\newcommand{\introthmname}{}
\newtheorem{introthminn}{\introthmname}

\newtheorem{thmx}{Theorem}
\newtheorem{corx}[thmx]{Corollary}

\newtheorem{propx}[thmx]{Proposition}

\usepackage{amssymb}
\usepackage{mathtools}
\usepackage{thmtools}
\usepackage{xfrac}
\usepackage{thm-restate}
\usepackage{enumerate}
\usepackage{hyperref}
\usepackage{xfrac}
\usepackage[capitalize]{cleveref}
\usepackage{verbatim}
\usepackage{amsmath}
\usepackage{amssymb}
\usepackage{mathabx}
\usepackage{graphicx}
\usepackage{amsfonts}
\usepackage{pb-diagram}
\usepackage{tikz-cd}
\usepackage{graphicx}
\usepackage{xcolor}
\usepackage{adjustbox}
\usepackage{mathtools}
\usepackage{todonotes}

\usepackage{mathrsfs}
\usepackage{bm}
\usepackage{enumitem}

\makeatletter
\newcommand{\tpitchfork}{%
  \vbox{
    \baselineskip\z@skip
    \lineskip-.52ex
    \lineskiplimit\maxdimen
    \m@th
    \ialign{##\crcr\hidewidth\smash{$-$}\hidewidth\crcr$\pitchfork$\crcr}
  }%
}

\DeclareMathOperator{\Homeo}{Homeo}

\newtheorem{theorem}{Theorem}
\newtheorem{definition}[theorem]{Definition}

\newtheorem{corollary}[theorem]{Corollary}
\newtheorem{lemma}[theorem]{Lemma}
\newtheorem{remark}[theorem]{Remark}
\newtheorem*{remark*}{Remark}
\newtheorem{observation}[theorem]{Observation}
\newtheorem{fact}[theorem]{Fact}

\newtheorem{proposition}[theorem]{Proposition}

\newcommand{\lu}[1]{{\color{magenta}#1}}

\newenvironment{subproof}[1][\proofname]{%
\begin{proof}[#1]%
	}{%
\end{proof}%
}

\newenvironment{subsubproof}[1][\proofname]{%
\begin{proof}[#1]%
}{%
\end{proof}%
}

\newcommand{\nn}[0]{\mathsf{N}}
\newcommand{\N}[0]{\mathbb{N}}

\newcommand{\Z}[0]{\mathbb{Z}}
\newcommand{\bx}[0]{\mathcal{B}(X)}

\newcommand{\V}[0]{\mathcal{V}}
\newcommand{\W}[0]{\mathcal{W}}

\DeclareRobustCommand{\loongrightarrow}{%
	\DOTSB\relbar\joinrel\relbar\joinrel\rightarrow
}

\DeclareMathOperator{\diam}{diam}
\DeclareMathOperator{\cdiam}{cdiam}
\DeclareMathOperator{\Fr}{Fr}
\DeclareMathOperator{\Int}{int}
\DeclareMathOperator{\supp}{supp}

\title[Group topologies on groups of bi-absolitely continuous homeomorphisms]{Group topologies on groups of bi-absolutely continuous homeomorphisms}
\author{J. de la Nuez Gonz\'alez }
\email{jnuezgonzalez@gmail.com}
\address{Korea Institute for Advanced Study (KIAS)}
\date{\today}
\thanks{Work supported by Samsung Science and Technology Foundation under Project Number SSTF-BA1301-51, by KIAS individual grant MG084001 and by the National Research Foundation of Korea under Mid-Career Research Program RS-2023-00278510.} 

\begin{document}

\newcommand{\h}[0]{\mathcal{H}^{+}_{bac}}  
\newcommand{\tac}[0]{\tau_{ac}}
\newcommand{\tc}[0]{\tau_{co}}
\newcommand{\idn}[0]{\mathcal{N}_{\tau}(\Id)}
\newcommand{\nm}[1]{\Vert#1 \rVert_{1}} 
\newcommand{\cf}[0]{\ind}
\renewcommand{\gg}[0]{\mathfrak{G}}
\renewcommand{\nn}[1]{\lVert #1 \rVert} 
\newcommand{\nac}[2]{\nn{#1}_{ac,#2}} 
\newcommand{\itg}[0]{L^{1}(X)}
\newcommand{\R}[0]{\mathbb{R}}
\newcommand{\no}[1]{\nn{#1}_{1}} 
\newcommand{\nd}[1]{\lVert #1 \rVert_{1,1}} 
\newcommand{\ms}[1]{\mu(#1)} 
\newcommand{\dm}[0]{\,d\mu}

\newcommand\ind{%
	\ifdefined\mathbbb%
	\mathbbb{1}%
	\else%
	\boldsymbol{\mathbb{1}}%
	\fi}

\numberwithin{theorem}{section}
\numberwithin{lemma}{section}
\numberwithin{corollary}{section}
\numberwithin{observation}{section}
\numberwithin{claim}{section}
\numberwithin{definition}{section}
\numberwithin{proposition}{section}
\numberwithin{fact}{section}
\numberwithin{remark}{section}

    \newtheorem{innercustomthm}{Theorem}
\newenvironment{customthm}[1]
{\renewcommand\theinnercustomthm{#1}\innercustomthm}
{\endinnercustomthm}

\newcommand{\hh}[0]{\mathcal{H}_{bac}}
\newcommand{\Id}[0]{\mathrm{Id}}

\begin{abstract}
	The group of homeomorphisms of the closed interval that are absolutely continuous and have an absolutely continuous inverse was shown by Solecki to admit a natural Polish group topology $\tau_{ac}$. We show that, under mild conditions on a compact space endowed with a finite Borel measure such a topology can be defined on the subgroup of the homeomorphism group consisting of those elements $g$ such that $g$ and $g^{-1}$ preserve the class of null sets. 
	
	 We use a probabilistic argument to show that in the case of a compact topological manifold equipped with an Oxtoby-Ulam measure, as well as in that of the Cantor space endowed with some natural Borel measures there is no group topology between $\tau_{ac}$ and the restriction $\tau_{co}$ of the compact-open topology. In fact, we show that any separable group topology strictly finer than $\tau_{co}$ must be also finer than $\tau_{ac}$. For one-dimensional manifolds we also show that $\tau_{co}$ and $\tau_{ac}$ are the only Hausdorff group topologies coarser than $\tau_{ac}$, and one can read our result as evidence for the non-existence of a good notion of regularity between continuity and absolute continuity.
	
	We also show that while Solecki's example is not Roelcke precompact, the group of bi-absolutely continuous homeomorphisms of the Cantor space endowed with the measure given by the Fr\"aiss\'e limit of the class of measured boolean algebras with rational probability measures is Roelcke precompact.  
\end{abstract}

\maketitle
This can be interpreted as showing the non-existence of a good notion of regularity between continuity and absolute continuity.

\section{Introduction}

\begin{definition}
	\label{d: absolute continuity}We say that $\phi:[0,1]\to\R$ is \emph{absolutely continuous} if  
  there exists some $\psi\in L^{1}([0,1])$ such that for all $t\in [0,1]$ we have 
		$\phi(t)-\phi(0)=\int_{[0,t]}\psi\,d\mu$.
\end{definition}

If $f$ is absolutely continuous, then the derivative of $f$ exists and coincides with $\psi$ almost everywhere. In \cite{solecki11999polish} Solecki considers the group
$$
 \h([0,1]):=\{\,g\in \mathcal{H}^{+}([0,1])\,\,|\,\,\text{$g$ and $g^{-1}$ are absolutely continuous}\,\},
$$
where $\mathcal{H}([0,1])$ denotes the homeomorphism group of a topological space $X$ ad $\mathcal{H}^{+}([0,1])$ the subgroup of those which are orientation preserving in case $[0,1]$ is a manifold, and he shows that it is a Polish group with the topology $\tau_{ac}$ induced by the metric $\rho(f,g):=\nn{f'-g'}_{1}$.
One can define the group $\hh(S^{1})$ analogously (see \cite{cohen2017polishability}, \cite{herndon2018absolute}) and prove Polishness in essentially the same way, although in this case an additional term controlling the displacement of a reference point is needed in the definition of $\rho$. 

Even though absolute continuity is particularly interesting in the $1$-manifold setting, in so far as it can be regarded as the first of an infinite family of intermediate regularity conditions defining some Polish subgroup of the $\mathcal{H}([0,1])$ (see \cite{cohen2017polishability}), we believe worth pointing out that Solecki's Polish group generalizes to a wide variety of other contexts, some of which will be explored here. 

Recall that given two non-negative measures $\mu,\nu$ on a measurable space $(X,\mathcal{A})$ we say that $\nu$ is absolutely continuous with respect to $\mu$, written $\nu\ll\mu$, if and only if $\nu(A)=0$ for any $A\in\mathcal{A}$ with $\mu(A)=0$. We are interested in the case where $\mu$ is a Borel measure, i.e., a measure on the Borel $\sigma$-algebra $\mathcal{B}(X)$ of $X$. 
For any such $\mu$ and any $g\in\mathcal{H}(X)$ there is a push-forward measure $g_{*}\mu$ given by $g_{*}\mu(A)=\mu(g^{-1}(A))$ for $A\in\mathcal{B}(X)$. We also write $g^{*}\mu:=(g^{-1})_{*}\mu$. 

\begin{definition}
	Let $X$ be a topological space and $\mu$ a Borel measure on $X$.  We let $\hh(X,\mu)$ be the subgroup of $\mathcal{H}(X)$ consisting of all those $g\in\mathcal{H}(X)$ such that $g_{*}\mu\ll\mu$ and $\mu\ll g_{*}\mu$ (equivalently, $g^{*}\mu\ll\mu$), that is, the subgroup of all $g\in\mathcal{H}(X)$ such that for $A\in\mathcal{B}(X)$ we have $\mu(A)=0$ if and only if $\mu(g(A))=0$.  
\end{definition}

Using a different terminology, $\hh(X,\mu)$ is the group of all homeomorphisms of $X$ that preserve the measure class of $\mu$. The condition $g_{*}(\mu)\ll\mu$ for a homeomorphisms in $g$ is satisfied precisely when $g$ is absolutely continuous in the sense of Banach, a terminology frequently found in older literature. We could thus refer to $\hh(X,\mu)$ as the set of bi-absolutely continuous homeomorphisms in the sense of Banach. \footnote{Note, however, that when $g$ is not a homeomorphism, continuity of $g$ and $g_{*}(\mu)\ll\mu$ are not sufficient conditions for absolute continuity. See \cite{rado2012continuous} or Chapter $14$ in \cite{bartle2001modern} for more details. } 


\begin{fact}[Radon-Nikodym, see \cite{klenke2013probability} Corollary 7.34]
	If $\mu,\nu$ are two $\sigma$-finite measures on a measurable space $(X,\mathcal{A})$ such that $\nu\ll\mu$, then there exists 
	a unique $\phi\in L^{1}(X,\mathcal{A},\mu)$ (that is, well-defined up to a $\mu$-null set), denoted by $\frac{d\nu}{d\mu}\geq 0$ such that $\nu(A)=\int_{A}\phi d\mu$ for all $A\in\mathcal{A}$, It follows that for all $\mathcal{A}$-measurable $h:X\to\R$ we have 
	$h\in L^{1}(X,\mathcal{A},\nu)$ if and only if $h\phi\in L^{1}(X,\mathcal{A},\mu)$ with equality
	$\int h\,d\nu=\int h\phi\,d\mu$ in that case. 
\end{fact}

If $g\in\hh(X,\mu)$, we write $\mathcal{D}(g):=\frac{dg^{*}\mu}{d\mu}$ (the generalized Jacobi determinant of $g$). We establish the following generalization of Solecki's result.
\newcommand{\extensionbody}[0]{
 	Let $X$ be a metrizable compact topological space and $\mu$ a finite Borel measure on $X$ such that $X$ has a basis of open sets with $\mu$-null frontier. Then the topology $\tau_{ac}$ on $\hh(X,\mu)$ generated by the restriction of the compact-open topology on $\mathcal{H}(X)$ and the topology induced by the distance $\rho_{ac}(g,h):=\nn{\mathcal{D}(g)-\mathcal{D}(h)}_{1}$ is a completely metrizable group topology. If in addition to this $X$ is separable, then $\tau_{ac}$ is Polish.
}

\begin{propx}
  \label{p: generalization}\extensionbody
\end{propx}

In case $X\in\{[0,1],S^{1}\}$ and $\mu$ is the Lebesgue measure we recover Solecki's Polish group (see \cite{bartle2001modern}, Theorem 14.15). Even in the original setting, our proof that $\tau_{ac}$ is a group topology diverges somewhat from that in \cite{solecki11999polish}. Our result applies, in particular, to compact manifolds of arbitrary dimension equipped with an Oxtoby-Ulam measure, as defined below.
\begin{definition}
	An Oxtoby-Ulam measure on a compact connected manifold $X$ is a finite Borel measure (necessarily regular) satisfying the following properties: 
	\begin{itemize}
		\item $\mu(U)>0$ for all non-empty open $U\subseteq X$,
		\item $\mu(\{x\})=0$ for all $x\in X$ 
		\item and $\mu(\partial X)=0$.
	\end{itemize}
  We will refer to $(X,\mu)$ for $X$ and $\mu$ as above as an OU pair.
\end{definition}
Another natural way of producing pairs $(X,\mu)$ as in the stement of the Proposition is the following. Suppose we are given an inverse system $\pi_{i}:X_{i+1}\to X_{i}$ of continuous epimorphisms between compact separable topological metric spaces, together with finite Borel measures $\mu_{i}$ on $X_{i}$ satisfying $(\phi_{i})_{*}(\mu_{i+1})=\mu_{i}$. Then by Theorem 2.2 in \cite{choksi1958inverse} the inverse limit $X:=\displaystyle{\lim_{\leftarrow}}\,X_{i}$  can be endowed with a limit Borel measure $\mu:=\displaystyle{\lim_{\leftarrow}}\,\mu_{i}$ such that for each of the natural projections $\pi_{i}:X\to X_{i}$ we have $(\pi_{i})_{*}(\mu)=\mu_{i}$.\footnote{Note that the inverse limit of a sequence of Borel measures does not exist in general, see \cite{choksi1958inverse} for more details.} 
If $X_{i}$ is generated by sets of $\mu_{i}$-null frontier for all $i$, it can be easily checked that the same property holds for $(X,\mu)$.

Two natural Borel measures on the Cantor space $\mathcal{C}$ are worth considering. Given some finite measure $\mu_{0}$ on $\{1,\dots r\}$ one can take the power $\mu_{0}^{\omega}$ measure on $\{1,\dots r\}^{\omega}\cong\mathcal{C}$ positive on singletons. We will refer to this as a generalized Bernoulli measure. 

Other than this, given any countable Field $K\subseteq\mathbb{R}$, one can consider the Fr\"aiss\'e limit $\mathbb{B}_{K}$ of the class of all the finite Boolean algebras (see, for instance, \cite{kechris2007turbulence} for $K=\mathbb{Q}$) equipped with probability measures with values on $K\cap[0,1]$. By Stone duality, one can identify $\mathbb{B}_{\mathbb{Q}}$ with the family of clopen sets of $\mathcal{C}$, with the measure on $\mathbb{B}_{K}$ extending to a unique (up to this identification) Borel measure on $\mathcal{C}$, which we will denote as $\mu_{K}$.  \\

\newcommand{\bodymainthm}[0]{
Let $(X,\mu)$ be either 
\begin{enumerate}[label=(\Roman*)]
	\item \label{ou case}an OU pair,
	\item \label{Cantor case}or the cantor space $\mathcal{C}=X$  equipped with either $\mu_{K}$ for some dense subfield $K\subseteq\mathbb{R}$ or a generalized Bernoulli measure. 
\end{enumerate}
Then there is no group topology $\tau$ on $\hh(X,\mu)$ with $\tau_{co}\subsetneq\tau\subsetneq\tau_{ac}$.
In fact, more is true: for any separable group topology $\tau$ on $\hh(X,\mu)$ if $\tau_{co}\subsetneq \tau$, then $\tau_{ac}\subseteq\tau$.
}

\subsection*{A gap between group topologies}  We will use the term $\tau_{co}$ to refer to the restriction of the compact-open topology to the group at hand. The main result of the paper is Theorem \ref{t: main} below, which establishes the existence of a gap between $\tau_{co}$ and $\tau_{ac}$ in a range of situations.   
\begin{thmx}
	\label{t: main}\bodymainthm
\end{thmx} 
A more general set of sufficient conditions for the first claim in the statement is given in Proposition \ref{p: first technical statement} and the reader interested only in this claim for the examples covered by Theorem \ref{t: main} may skip the discussion in Section \ref{s: separability}. 

In the case of a one-dimensional manifold $X$, one can read the result above as stating that there is no notion of regularity (for homeomorphisms) sitting between between continuity and absolute continuity in the same way absolute continuity sits between $C^{1}$-regularity and continuity.

In the one dimensional case we actually know that 
$\tau_{co}$ is the coarsest (not just minimal) Hausdorff group topology on $\hh(X,\mu)$ by the arguments in \cite{gartside2003autohomeomorphism}. Therefore, we can conclude:
\begin{corx}
	Let $(X,\mu)$ be an OU pair with $X\in\{[0,1],S^{1}\}$ and $\mu$ be as above. Then $\tau_{ac}$ and $\tau_{co}$ are the Hausdorff group topologies on $\hh(X,\mu)$ coarser than $\tau_{co}$. 
\end{corx}

  \subsection*{Roelcke precompactness} We say that a topological group $(G,\tau)$ is Roelcke precompact if it is totally bounded with respect to the Roelcke uniformity or, in more concrete terms, if for any neighborhood $\V$ of the identity in $G$ there are finitely many elements $\{g_{1},\dots g_{r}\}$ such that $G=\bigcup_{l=1}^{r}\V g_{l}\V$.
Roelcke precompactness is the topological correlate of the model theoretical notion of $\omega$-categoricity, which explains some of its importance. 

In \cite{rosendal2021coarse} Rosendal develops a theory that endows every topological group in some large family that contains all Polish groups with a canonical metrizable coarse structure. Roelcke precompactness implies boundedness in the sense of said coarse structure. It is shown in \cite{herndon2018absolute} that the group $(\hh(X,\mu),\tau_{ac})$ where $X\in\{[0,1],S^{1}\}$ and $\mu$ is the standard measure, is bounded. One might wonder whether this could just be a consequence of $\hh(X,\mu)$ being Roelcke precompact. However, we show the following.

\newcommand{\mainthmbodyy}[0]{
Let $(X,\mu)$ be a $1$-dimensional OU pair. Then $\hh(X,\mu)$ is not Roelcke precompact. More precisely, any sequence $(g_{n})_{n\in\N}\subset\hh(X,\mu)$ such that the following conditions are satisfied: 
\begin{itemize}
	\item $(g_{n})_{n\in\N}$ converges to the identity in $\tau_{co}$,
	\item and $(g_{n})_{n\in\N}$ does not converge to $\Id$ in $\tau_{ac}$   
\end{itemize} 
is not totally bounded for the Roelcke uniformity. In particular $\hh(X,\mu)$ is not Roelcke precompact.
} 
\begin{thmx}
	\label{t: main2}\mainthmbodyy
\end{thmx} 

\begin{remark}
	In the case where $(X,\mu)$ is an OU pair with $dim(X)\geq 2$, the non Roelcke precompactness of $\hh(X,\mu)$ follows from that of $(\mathcal{H}(X),\tau_{co})$ (see \cite{rosendal2013global}) together with Corollary \ref{c: density}.
\end{remark}

One might wonder whether it might be possible to abstract the proof of Theorem \ref{t: main2} so as to show that the conclusion holds for any sufficiently non-trivial instance of Proposition \ref{p: generalization}. However, some degree of uniform control of the topology of $X$ over the measure seems to be needed, as evidenced by our last result, which plays out in the opposite direction. 

\newcommand{\mainthmbodyyy}[0]{
Let $K$ a countable subfield of $\mathbb{R}$. Then $(\hh(\mathcal{C},\mu_{K}),\tau_{ac})$ is Roelcke precompact.
} 
\begin{thmx}
	\label{t: main3}\mainthmbodyyy
\end{thmx} 

\subsection*{Structure of the paper}

\cref{s:polish} is dedicated to the proof of Proposition \ref{p: generalization}, while introducing useful notation and results used later on. \cref{s:OUpairs} focuses on the properties of $\mathcal{H}_{bac}(X,\mu)$ in the case in which $(X,\mu)$ is an OU pair. It collects technical statements generalizing properties of the group of homeomorphisms of a compact manifold, such as a fragmentation property relative to $\tau_{co}$. \cref{s:preliminaries} collects a handful of miscellaneous results and definitions that will appear in the two heavier sections after it. 

Sections \ref{s: well distributed elements} and \ref{s:dense large elements} both unfold in the setting in which we are given $(X,\mu)$ as in Theorem \ref{t: main} and some group topology $\tau$ on $\hh(X,\mu)$ with $\tau_{co}\subseteq\tau$, $\tau_{ac}\nsubseteq\tau$. The main result of \cref{s: well distributed elements} is \cref{c: uniform derivative}, which states that in this setting, given any finite collection of disjoint sets $\{A_{i}\}_{i=1}^{k}$ in $X$ every neighborhood of the identity in $\tau$ must contain some element $g$ that is uniformly bounded away from the identity in each $A_{i}$, in the sense of $\rho_{ac}$. \cref{s:dense large elements} uses this to derive \cref{c: shrinking sets}, which asserts that in the above setting, every neighborhood of the identity contains some element that maps an open set $U$ with $\mu(U)$ arbitrarily close to $\mu(X)$ to a set of arbitrarily small measure. 

\cref{s: separability} considers a group topology $\tau$ in $\hh(X,\mu)$ of which it is only assumed to be separable in the case in which $(X,\mu)$ is an OU pair. Adapting well known arguments found in the literature on automatic continuity of homeomorphism groups yields conclusions which fall short of automatic continuity for $\hh(X,\mu)$, but which enable the final stronger claim in Theorem \ref{t: main}. The theorem itself is proved in \cref{s:proof of b}, while \cref{s:roelcke} is devoted to the proof of Theorem \ref{t: main2} and Theorem \ref{t: main3}. 

\section{A family of Polish groups}
\label{s:polish}

Let $X$ be a compact metric space and $\mu$ a finite Borel measure on $X$. We will assume throughout the paper that in such situation some compatible metric $d$ on $X$ has been fixed. 
For readability, in the presence of integrals we will often use $\circ$ to denote the composition of elements in $\mathcal{H}(X)$ and write $\Id$ for the identity map. We write $\ind_{A}$ for the indicator function of a set $A\subseteq X$. 
For convenience, we will often write $1$ instead of $\ind_{X}$.  
The following two Observations are well-known and easy to show (see \cite{solecki11999polish}). We include a proof of the next one for completeness. 
\begin{observation}
	\label{o: change of variable} For any $g\in\hh(X,\mu)$ and any $\phi\in L^{1}(X,\mu)$ we have 
	$$
	\int (\phi\circ g)\mathcal{D}(g)\,d\mu=\int \phi\,d\mu
	$$
\end{observation}
 \begin{proof}
 	$$\int (\phi\circ g)\mathcal{D}(g)\,d\mu=\int \phi\circ g\,dg^{*}\mu =\int\phi\circ g\circ g^{-1}\,d\mu=\int\phi\,d\mu $$
 \end{proof}
Using this and the fact that $\mathcal{D}(g)$ is well-defined up to a set of measure zero, one also checks: 
\begin{observation}
	\label{o: product rule} For any $g_{1},g_{2}\in\hh(X,\mu)$ we have 
	\begin{equation*}
		\mathcal{D}(g_{1}\circ g_{2})=(\mathcal{D}(g_{1})\circ g_{2})\mathcal{D}(g_{2}) \quad\quad \mathcal{D}(g_{1}^{-1})=\frac{1}{\mathcal{D}(g_{1})\circ g_{1}^{-1}},\quad\quad\text{a.e.}
	\end{equation*}
\end{observation}

Recall that given $g_{1},g_{2}\in\hh(X,\mu)$ we write $$\rho_{ac}(g_{1},g_{2}):=\nn{\mathcal{D}(g_{1})-\mathcal{D}(g_{2})}_{1} \quad\quad \nn{g_{1}}_{ac}:=\rho_{ac}(g_{1},\Id).$$
Clearly $\rho_{ac}$ is a pseudo-metric on $\gg$. For $A\in\mathcal{B}(X)$,  $f\in L^{1}(X,\mu)$ and $g\in\gg$ we also write
$$\nn{f}_{1,A}:=\nn{\ind_{A}f}_{1} \quad\quad \nn{g}_{ac,A}:=\nn{\mathcal{D}(g)-1}_{1,A}.$$

Write $\tau_{ac}:=\tau_{ac}(X,\mu)$ for the topology generated by the restriction to $\hh(X,\mu)$ of the compact-open (uniform convergence) topology on $\mathcal{H}(X)$, which we shall denote henceforth simply as $\tau_{co}$, and the topology induced by the pseudo-metric $\rho_{ac}$ above.

Fix a compatible metric $d$ on $X$ (necessarily complete).
 \begin{definition}
	 \label{d: neighborhoods}Given $\delta>0$ we write
	 $$
	  \V_{\delta}:=\{g\in\gg\,|\,\sup_{x\in X}d(g(x),x)\leq\delta\},\quad\quad\V'_{\delta}:=\{g\in\gg\,|\,\nn{g}_{ac}\leq\delta\},\quad\quad \W_{\delta}:=\V_{\delta}\cap\V'_{\delta},
	 $$
 \end{definition}
 so that a system of neighborhoods of $\Id$ in $\tau_{ac}$ is given by $\{\W_{\delta}\}_{\delta>0}$.

 The following observation is well-known (the first statement is a formal generalization of a result in \cite{solecki11999polish}), but we include a proof for completeness.
 \begin{observation}
  \label{o: right invariance} For any $h_{1},h_{2},g\in\gg$ we have 
  $
  \rho_{ac}(h_{1},h_{2})=\rho_{ac}(h_{1}\circ g,h_{2}\circ g)
  $. Therefore, the following equalities hold:
   \begin{equation*}
   	\nn{h_{1}}_{ac}=\nn{h_{1}^{-1}}_{ac},\quad \quad \nn{h_{1}\circ h_{2}}_{ac}\leq\nn{h_{1}}_{ac}+\nn{h_{2}}_{ac}.
   \end{equation*}
 \end{observation}
 \begin{proof}
	The first equality is an easy consequence of Observations \ref{o: change of variable} and \ref{o: product rule},
	while the second equality follows suit by taking $h_{2}=\Id,g=h_{1}^{-1}$. For the final inequality:
	$$
	\nn{h_{1}\circ h_{2}}_{ac}=\rho_{ac}(h_{1},h_{2}^{-1})\leq \rho_{ac}(h_{1},1)+\rho_{ac}(h_{2}^{-1},1)=\nn{h_{1}}_{ac}+\nn{h_{2}}_{ac},
	$$
	where the second equality above follows from the equality $\nn{h_{2}^{-1}}_{ac}=\nn{h_{2}}_{ac}$. 
 \end{proof}

\begin{definition}
	\label{d: continuity modulus}Given $f\in C^{0}(X)$ we write 
	$$
	\omega_{f}:\R_{>0}\to\R_{\geq 0} \quad\quad \omega_{f}(\delta)=\sup\{\,|f(x)-f(y)|\,|\,d(x,y)\leq\delta\,\}
	$$
\end{definition}

It is worth highlighting the following observation, which will be used throughout the paper:
\begin{observation}
	\label{o: approximation use} Let $g\in\gg$ and $A\in\mathcal{B}(X)$. If $\phi\in C^{0}(X)$ satisfies 
	$\nn{\phi-\mathcal{D}(g)}_{1,A}\leq\epsilon$ and $\omega_{\phi}(\delta)\leq\epsilon'$ then for any $h\in\V_{\delta}$ we have: 
	$$
	\nn{\mathcal{D}(g\circ h)-\phi \mathcal{D}(h)}_{1,A}\leq\epsilon+\epsilon'\mu(h(A)).
	$$ 
	If $\omega_{\phi}(2\delta)\leq\epsilon'$ and $\psi:=\sum_{i=1}^{m}\lambda_{i}\ind_{A_{i}}$ for some partition $\{A_{i}\}_{i=1}^{m}$ of $X$ into Borel sets of diameter at most $\delta$ and $\lambda_{i}\in[\inf_{A_{i}}\phi,\sup_{A_{i}}\phi]$, then the inequality above holds with $\psi$ in place of $\phi$. 
\end{observation}
\begin{proof}
	Using the composition rule and the triangular inequality for the first inequality and Observation \ref{o: change of variable} in the equality after that we obtain 
	\begin{align*}
		\nn{\mathcal{D}(g\circ h)-\phi \mathcal{D}(h)}_{1,A}\leq& \nn{((\mathcal{D}(g)-\phi)\circ h)\mathcal{D}(h)}_{1,A}+\nn{(\phi\circ h-\phi)\mathcal{D}(h)}_{1,A}=\\
		&\nn{\mathcal{D}(g)-\phi}_{1,A}+\nn{(\phi\circ h-\phi)\mathcal{D}(h)}_{1,A}\leq\\
		&\epsilon+\nn{\phi\circ h-\phi}_{\infty}\nn{\mathcal{D}(h)}_{1,A}\leq\epsilon+\epsilon'\mu(h(A)).	
	\end{align*}
	The second claim follows by a similar argument, since the assumption that $\omega_{\phi}(2\delta)<\epsilon'$ implies that 
	$\nn{\phi\circ h-\psi}_{\infty}\leq\epsilon'$. 
\end{proof}

\begin{lemma}
	\label{l: continuity on the right}For any $g\in\gg$ and $\epsilon>0$ there is $\delta=\delta(\epsilon,g)$ such that for any 
	$h\in\W_{\delta}$ we have $\nn{h^{g^{-1}}}_{ac}=\rho_{ac}(g,g\circ h)\leq\epsilon$.  
\end{lemma}
\begin{proof}
	Choose some $\phi\in C^{0}(X)$ such that $\nn{\phi-\mathcal{D}(g)}_{1}<\frac{\epsilon}{4}$ and $\delta>0$ such that  
	$$\delta\leq\min\left\{\,\omega_{\phi}(\frac{\epsilon}{4}),\,\frac{\epsilon}{4\nn{\phi}_{\infty}}\,\right\}$$.
	Then using \cref{o: approximation use} in the second inequality below, we get 
	\begin{align*}
	   \nn{\mathcal{D}(g)-(\mathcal{D}(g)\circ h)\mathcal{D}(h)}_{1}\leq&\nn{\phi-\phi \mathcal{D}(h)}_{1}+\nn{\phi-\mathcal{D}(g)}+\nn{((\phi-\mathcal{D}(g))\circ h)\mathcal{D}(h)}_{1}
	   \leq\\
	   &\nn{\phi-\phi \mathcal{D}(h)}_{1}+\frac{\epsilon}{4}+(\frac{\epsilon}{4}+\frac{\epsilon}{4})\leq\nn{\phi}_{\infty}\nn{\mathcal{D}(h)-1}_{1}+\frac{3\epsilon}{4}\leq\epsilon,
	\end{align*}
  as needed. 
\end{proof}
  
\begin{corollary}
	Under the assumptions on $X$ and $\mu$ made at the beginning of the section, $\tau_{ac}$ is a group topology on $\gg$. 
\end{corollary}
\begin{proof}
  \cref{l: continuity on the right}, together with the fact that $\tau_{co}$ is a group topology implies that the unique left-invariant topology on $\gg$ for which $\{\W_{\delta}\}_{\delta>0}$ is a system of neighborhoods of the identity is in fact a group topology. It has to coincide with $\tau_{ac}$ by the same Lemma. 
 \end{proof}

It remains to show that $(\gg,\tau_{ac})$ is Polish.
Consider now the metric $\rho:\gg \times\gg \to\R_{\geq 0}$ on $\gg$ given by: 
$$
 \rho(g_{1},g_{2})=\sup_{x\in X}d(g_{1}(x),g_{2}(x))+\sup_{x\in X}d(g_{1}^{-1}(x),g_{2}^{-1}(x))+\rho_{ac}(g_{1},g_{2})+\rho_{ac}(g_{1}^{-1},g_{2}^{-1}).
$$ 

This induces the topology $\tau_{ac}$ on $\gg$ by Observation \ref{o: right invariance}, together with the fact that the metric given by the sum of the first two terms above is a complete metric for the compact-open topology.

We note the following fact:  
\begin{fact}
	\label{f: continuity} Suppose we are given two finite measures $\mu,\nu$ on a measurable space $(X,\mathcal{A})$ for which
	 $\nu\ll\mu$. Then for any $\epsilon>0$ there exists some $\delta>0$ such that $\mu(g(A))<\epsilon$ for all 
	$A\in\mathcal{A}$ with $\mu(A)<\delta$.  
\end{fact}
Given $g\in\gg$ and $\delta>0$ we write
$$
\Omega_{g}(\delta):=\sup\,\{\mu(g(A))\,|\,A\in\mathcal{B}(X),\,\,\mu(A)\leq\delta) \}
$$

\begin{lemma}
	\label{l: uniform continuity in measure} Let $\underline{g}=(g_{n})_{n\geq 0}$ be a sequence in $\gg$ that is Cauchy with respect to $\rho$. Then the function 
	$\Omega_{\underline{g}}=\sup_{n\geq 0}\Omega_{g_{n}}$ satisfies $\displaystyle{\lim_{\epsilon\to 0}}\,\Omega_{\underline{g}}(\epsilon)=0$. 
\end{lemma}
\begin{proof}
	For a fixed $\epsilon>0$, take $n_{0}$ such that for all $n\geq n_{0}$ we have $\nn{\rho_{ac}(g_{n},g_{n_{0}})}\leq\frac{\epsilon}{2}$.
	Then $\sup_{n\geq 0}\mu(g_{n}(A))\leq\epsilon$ for any $A\in\mathcal{B}(X)$ such that $\Omega_{g_{k}}(\mu(A))\leq\frac{\epsilon}{2}$ for all $0\leq k\leq n_{0}$. 
\end{proof}

\begin{lemma}
	\label{l: completeness} $\rho$ is a complete metric on $\gg$.
\end{lemma}
\begin{proof}
	Let $\underline{g}=(g_{n})_{n\geq 0}$ be a Cauchy sequence in $\gg$ with respect to $\rho$. Then $(g_{n})_{n\geq 0}$ must converge in $\tau_{co}$ to some $g_{\infty}\in\mathcal{H}(X)$, by completeness of the metric $\rho_{\infty}(g,h)=\sup_{x\in X}d(g_{1}(x),g_{2}(x))+\sup_{x\in X}d(g^{-1}_{1}(x),g^{-1}_{2}(x))$ on $\mathcal{H}(X)$. We claim that $g_{\infty}\in\hh(X,\mu)$.   
  Note that $(\mathcal{D}(g_{n}))_{n\geq 0}$ must be Cauchy in $L^{1}(X,\mu)$ and thus 
  must converge in $L^{1}(X,\mu)$ to some function $\Psi\geq 0$, which implies for any $A\in\bx$ we have: 
   \begin{equation}
   	\label{e: convergence measure}
   	 \mu(g_{n}(A))=\int_{A}\mathcal{D}(g_{n})\,d\mu\underset{n\to\infty}{\loongrightarrow}\int_{A}\Psi \,d\mu.
   \end{equation}
  
 	Since $g_{n}\to g_{\infty}$ in $\tau_{co}$ (i.e., the topology of uniform convergence), for any open set $U\subseteq X$ we must have 
	$g_{\infty}(U)\subseteq\displaystyle{\bigcup_{n\geq 0}\bigcap_{k\geq n}} g_{k}(U)$. The continuity of $\mu$ then allows us to conclude
	\begin{equation}
		\label{e: upper continuity on measure}
		\mu(g_{\infty}(U))\leq\lim_{n\to \infty}g_{n}(U) 
	\end{equation} 
	and so, in particular, $\mu(g_{\infty}(U))\leq\Omega_{\underline{g}}(\mu(U))$ for any such $U$, which implies that $\mu(g_{\infty}(A))=0$ for any $A\in\mathcal{B}(X)$ with $\mu(A)=0$, by Lemma \ref{l: uniform continuity in measure} and the regularity of $\mu$. A similar observation can be made for the inverse $g_{\infty}^{-1}$ and the (Cauchy) sequence $(g_{n}^{-1})_{n\geq 0}$ and so we conclude that $g_{\infty}\in \hh(X,\mu)$.  
	
	If now $U\subseteq X$ is an open set with
	$\mu(Fr(U))=0$ then $\mu(Fr(g_{\infty}(U)))=0$ and thus by (\ref{e: upper continuity on measure}) we have
	\begin{equation*}
		\mu(X)=\mu(g_{\infty}(U))+\mu(g_{\infty}(X\setminus \bar{U}))\leq\lim_{n\to\infty}[\mu(g_{n}(U))+\mu(g_{n}(X\setminus\bar{U}))]=\mu(X),
	\end{equation*}
  which implies 
  $$\mu(g_{\infty}(U))=\lim_{n\to\infty}\mu(g_{n}(U))=\int_{U}\Psi\,d\mu$$
   for all such $U$. Note that any finite union of open sets with $\mu$-null frontier also has $\mu$-null frontier.
   Since by assumption the topology of $X$ is generated by open sets with $\mu$-null frontier, it follows that every closed (i.e compact) set $F\subseteq X$ is the intersection of a descending sequence of open sets $U_{0}\supset U_{1}\dots$ with $\mu$-null frontier. To see this, observe that for all $n$ the set $F$ will be covered by a finite family of open sets of diameter at most $\frac{1}{n}$ with $\mu$-null frontier and take as $U_{n}$ their union, which also has $\omega$-null frontier. Using the dominated convergence theorem, we now have
   $$
     \mu(g_{\infty}(F))=\lim_{n}\mu(g_{\infty}(U_{n}))=\lim_{n}\int_{U_{n}}\Psi\dm=\int_{F}\Psi\dm.
   $$  
   The regularity of the measures $g_{\infty}^{*}\mu$ and $A\mapsto\int_{A}\Psi\dm$ now easily shows that $\mu(g_{\infty}(A))=\int_{A}\Psi\dm$ for all $A\in\mathcal{B}(X)$ and thus that $\mathcal{D}(g_{\infty})=\Psi$ and $\mathcal{D}(g_{n})$ converge to $\mathcal{D}(g_{\infty})$ in $L_{1}$. Arguing similarly for $g_{\infty}^{-1}$ and $(g_{n}^{-1})_{n\geq 0}$,
  we conclude that $\lim_{n\to\infty}\rho(g_{n},g_{\infty})=0$.
\end{proof}

To conclude the proof of Proposition \ref{p: generalization} we note that the separability of $X$ implies that of $\hh(X,\mu)$. 
\begin{observation}
	\label{o: separability}If $X$ is separable, then the topology $\tau_{ac}$ is separable.
\end{observation}
\begin{proof} 
   Let $\{\phi_{n}\}_{n\geq 1}$ some dense subset of $L^{1}(X,\mu)$ and $\{h_{n}\}_{n\geq 1}$ 
	countable dense subset of $\mathcal{H}(X)$ for the compact-open topology. For every pair of integers $n,m,N\geq 1$ let 
	$g_{n,m,N}\in\hh(X,\mu)$ be either an element such that $\nn{\mathcal{D}(g_{n,m,N})-\phi_{n}}_{1}\leq\frac{1}{N}$ and 
	$\sup_{x\in X}d(g_{n,m,N}(x),h_{n}(x))\leq\frac{1}{N}$ in the case that such an element exists and arbitrary otherwise. It is easy to check that $\{g_{n,m,N}\}_{n,m,N\geq 1}$ is dense in $\tau_{ac}$ and we leave it to the reader. 
\end{proof}

\section{Some properties of $\hh(X,\mu)$ for an Oxtoby-Ulam pair $(X,\mu)$.}
  \label{s:OUpairs}
   \newcommand{\cs}[0]{\widebar{supp}}
   \newcommand{\mh}[0]{\mathfrak{h}}
   
   Fix an OU pair $(X,\mu)$. In what follows, we will write $\gg$ for $\hh(X,\mu)$.
   By a \emph{nicely embedded ball} in $X$ we mean a homeomorphic image in $int(X)$ of the standard closed $m$-ball 
   $\bar{B}^{m}(1)$ that can be extended to a homeomorphic embedding of the open ball $B^{m}(2)$. By a \emph{boundary half-ball} we mean the image of $\bar{B}^{m}(1)\cap\R^{m-1}\times\R_{\geq 0}$ by a homeomorphic embedding $h:\bar{B}^{m}(1)\cap\R^{m-1}\times\R_{\geq 0}$ which extends to a homeomorphic embedding $\tilde{h}:B^{m}(2)\cap\R^{m-1}\times\R_{\geq 0}\to X$ with the property that $\tilde{h}^{-1}(\partial X)=B^{m}(2)\cap\{x_{m}=0\}$. By a \emph{$\mu$-adapted ball (resp. boundary half ball)} we mean a nicely embedded ball (resp. boundary half ball) whose frontier is $\mu$-null. 
   
   We denote by $\mathcal{H}(X,\partial X)$ the group of homeomorphisms of $X$ fixing $\partial X$. 
   Note that this is still a Polish group when endowed with the compact-open topology.  
   The following classical result is due to Oxtoby and Ulam \cite{oxtoby1941measure}.
   \begin{fact}
   	\label{f: OU uniqueness}Let $X$ be a compact manifold and $\mu,\mu'$ two OU measures on $X$ with $\mu(X)=\mu'(X)$. Then there exists some $h\in\mathcal{H}(X,\partial X)$ such that $h_{*}\mu=\mu'$. 
   \end{fact}
 
   Given $g\in\mathcal{H}(X)$, write $\supp(g):=\{x\in X\,|\,g(x)\neq x\}$ and denote its closure by $\cs(g)$.
   Given an open set $U\subseteq X$ and $\mathfrak{h}\leq\mathcal{H}(X)$ we write $$\mathfrak{h}_{U}:=\{g\in\mathfrak{h}\,|\,\supp(g)\subseteq U\},\quad\quad\quad \mathfrak{h}_{U}^{c}:=\{g\in\mathfrak{h}\,|\,\cs(g)\subseteq U\} .$$

   Given $A\in\mathcal{B}(X)$ we say that $h\in\mathcal{H}(X)$ is measure-class preserving on $A$ if for all $C\in\mathcal{B}(X)$, $C\subseteq A$ we have that $\mu(C)=0$ if and only $\mu(h(C))=0$. 
   The following is obvious from $\sigma$-additivity:  
  \begin{observation}
  	\label{o: pasting}Let $F_{k}\in\bx$, $k\geq 0$
  	then $g\in\mathcal{H}(X)$ is measure-class preserving on $\bigcup_{k\geq 0}F_{k}$ if it is on $F_{k}$ for all $k\geq 0$. 
  \end{observation}

   \begin{observation}
   	\label{o: ball correction for any}Let $h\in\mathcal{H}(X)$ and assume that $Y$ is a closed subspace of $X$ homeomorphic to a compact $m$-manifold such that $\mu(\partial Y)=0$. Then there exists some $g\in\mathcal{H}(X)$ supported on $Y$ such that $hg$ is measure-class preserving on $\bigcup_{j\in\mathcal{F}}Y_{j}$. Moreover, we may choose $g$ so that
     $$
     \nn{hg}_{ac,Y}=|\mu(Y)-\mu(h(Y))|.
     $$    
   \end{observation}
   \begin{proof}
 Apply \cref{f: OU uniqueness} to the two OU measures $\nu_{1}:=\frac{1}{\mu(h(Y))}h^{*}(\mu_{\restriction h(Y)})$ and $\nu_{2}:=\frac{1}{\mu(Y)}\mu_{\restriction Y}$ on the compact manifold $Y$. This yields an element $\gamma\in\mathcal{H}(Y,\partial Y)$ such that $\gamma^{*}(\nu_{1})=\nu_{2}$. We then have 
   	\begin{equation*}
   		(h\circ\gamma)^{*}(\mu_{\restriction h(Y)})=\mu(h(Y))\gamma^{*}(\nu_{1})=\mu(h(Y))\nu_{2}=\frac{\mu(h(Y))}{\mu(Y)}\mu_{\restriction Y}.
   	\end{equation*}
   	Therefore, if we extend $\gamma$ by the identity outside $Y$ to an element of $\mathcal{H}(X)$, then 
   	$hg$ is measure-class preserving on $Y$. Thus, $\mathcal{D}(hg)$ is well-defined on $Y$ and equal to the constant function $\frac{\mu(h(Y))}{\mu(Y)}$, so that $\nn{hg}_{ac,Y}$ is also well-defined and equal to $|\mu(Y)-\mu(h(Y))|$. 
   \end{proof}
    
   The following is most likely known, but we have failed to find an appropriate reference. 
   
   \begin{corollary}
   	\label{o: ball transitivity}Let $\{B_{i}\}_{i=1}^{r}$, $\{B'_{i}\}_{i=1}^{r}$ be two collections of disjoint $\mu$-adapted balls and boundary half-balls $X$ with the property that $B_{i}$ is a boundary half-ball if and only if $B'_{i}$ is, in which case they both intersect the same boundary component. If $m=2$, then we also assume that the intervals in which $B_{i}$ and $B'_{i}$ touch the boundary components of $X$ are ordered identically for the two collections. Then there exists $h\in\gg$ such that $h(B_{i})=B'_{i}$ for all $1\leq i\leq r$.  
   \end{corollary}
   \begin{proof}
   	It follows from the work of Kirby \cite{kirby2010stable} and Quinn \cite{quinn1982ends} that there exists $h_{0}\in\mathcal{H}(X)$ 
   	with $h_{0}(B_{i})=B'_{i}$ for all $1\leq i\leq r$. The result then follows by applying \cref{o: ball correction for any} to each of the members of the covering of $X$ given by $\{B_{i}\}_{i=1}^{r}\cup\{\overline{X\setminus\bigcup_{i=1}^{r}B_{i}}\}$. For each member $A$ of this covering, we thus obtain an element $g_{A}$ supported on $A$ with the property that $h_{0}\circ g_{A}$ is measure-class preserving on $A$. The product of all $g_{A}$ (which commute) is an element $g$ with the property that $h:=h_{0}g\in\gg$ and acts on the balls as needed.  
    \end{proof}

   	  For the following, see Claim 3.12 in \cite{mann2016automatic}:
   	\begin{fact}
   		\label{f: colorable coverings} Given a compact manifold $X$ there is some $N=N(X)>0$ such that for any $\epsilon>0$ we can write $X=\bigcup_{l=1}^{N}V_{l}$, where $V_{l}$ is a disjoint union of interiors of embedded balls and boundary balls of diameter less than $\epsilon$.
   	\end{fact}

  \begin{lemma}
  	 \label{l: nice colorable coverings}In the conclusion of \ref{f: colorable coverings} we may assume that all balls and boundary half-balls involved are $\mu$-adapted. 
  \end{lemma}
   \begin{proof}
  By using the fact that the collection $\{\mathring{B}_{i,\lambda}\,|\,1\leq j\leq N,\lambda\in\Lambda_{j}\}$ is locally finite and considering the Lebesgue number of the restriction of the covering above to different compact sets, it is not hard to find a function $\mathbf{\delta}:X\to\R_{>0}$ bounded away from zero on compact sets and such that if one replaces each of the balls or boundary half-ball $B$ in the collections by a smaller one $B'$ of the same type and with the property that the Hausdorff distance between $\Fr(B)$ and $\Fr(B')$ is at most $\inf_{x\in B}\delta(x)$, then the resulting new collections will cover $X$. That it is possible to do so in a way that the resulting (embedded) balls are $\mu$-adapted follows from Fact \ref{f: OU uniqueness}.  
   \end{proof}
 
   \newcommand{\tv}[0]{\tilde{\V}}
	 \renewcommand{\H}[0]{\mathfrak{h}}
	  \newcommand{\Nb}[0]{\mathrm{N}}
   
   \begin{lemma}
   	\label{l: ac correction} Let $g\in\mathcal{H}(X)$ and let $U$ be an open subset of $X$, $F\subseteq U$ a compact set, and $\epsilon>0$. 
   	
   	Then there exists some $h\in\mathcal{H}(X,\partial X)$ satisfying:
   	\begin{itemize}
   		\item $F\subseteq\overline{supp}(h)\subseteq\supp(h)\subseteq U$, 
   		\item $h\in\tv_{\epsilon}$,
   			\item $gh$ is measure-class preserving on $\overline{\supp}(h)$.
   	\end{itemize}
   \end{lemma}
     \begin{proof}
   	Let $N(X)$ be as in \cref{f: colorable coverings}. Using \cref{l: nice colorable coverings} construct $N$ finite families $\{B_{i}^{l}\}_{i\in\mathcal{J}_{l}}$ of $\mu$-adapted balls and boundary half balls such that $\diam(B_{i}^{l})<\frac{\epsilon}{N}$ for all $1\leq l\leq N$ and $i\in\mathcal{J}_{l}$ and the union of the interiors of all the $B_{i}^{l}$ covers $X$. 
   	
   	Let $C_{1},\dots C_{m}$ be an enumeration of all sets of the form $B^{l}_{i}$ that intersect $F$. Clearly, $F\subseteq\bigcup_{j=1}^{m}\mathring{C}_{j}$ and, by taking $\epsilon$ to be small enough, we may assume that $\bigcup_{j=1}^{m}C_{j}\subset U$. 
   	
   	Using \cref{o: ball correction for any} we can find $h_{1}\in\mathcal{H}(X,\partial X)$ with $\supp(h_{1})\subseteq C_{1}$ such that 
   	$gh_{1}$ is measure-class preserving on $C_{1}$. Applying \ref{o: ball correction for any} once more (to the homeomorphism $g'=gh_{1}$) yields $h_{2}\in\mathcal{H}(X,\partial X)$ such that $gh_{1}h_{2}$ is measure-class preserving on $C_{2}$. Since $(gh_{1})_{\restriction C_{1}\setminus C_{2}}=(gh_{1}h_{2})_{\restriction C_{1}\setminus C_{2}}$, $gh_{1}h_{2}$ is also measure-class preserving on $C_{1}\setminus C_{2}$. By \cref{o: pasting} $gh_{1}h_{2}$ is thus measure-class preserving on $C_{1}\cup C_{2}$. Continuing in this fashion, we obtain $h_{1},\dots h_{m}\in\mathcal{H}(X,\partial X)$ with $h_{j}$ supported on $C_{j}$ and such that if we write $h:=h_{1}h_{2}\cdots h_{m}$, then the composition $gh$ is measure-class preserving on $\bigcup_{j=1}^{m}C_{j}$, which is a neighborhood of $F$ contained in $U$. 
   	
   	Since $\diam(C_{j})<\frac{\epsilon}{N}$ and one can still separate the family $\{C_{j}\}_{j=1}^{m}$ into $N$ distinct families of disjoint sets, it is not hard to see that $h\in\tv_{\epsilon}$ and so we are done.  
   \end{proof}
   
   \begin{corollary}
   	\label{c: density} For any OU pair $(X,\mu)$ the group $\hh(X,\mu)$ is dense in $\mathcal{H}(X)$ with respect to $\tau_{co}$. 
   \end{corollary}

 \subsection*{Local fragmentation for the restriction of the compact-open topology}

    For $\delta>0$ write $\tv_{\delta}:=\{g\in\mathcal{H}(X)\,|\,\forall x\in X\,d(x,g(x))\leq\delta\}$.  
    We write $\Nb_{\epsilon}(A)$ for the metric neighborhood $\Nb_{\epsilon}(A):=\{x\in X\,|\,d(x,A)<\epsilon\}$. 
   
   \begin{definition}
   	\label{d: fragmentation} We say that a subgroup $\mathfrak{h}\leq\mathcal{H}(X)$ satisfies the weak local fragmentation property for $\tau_{co}$ if for any finite collection of open sets $U_{1},\dots U_{r}$ and any open $V\subseteq X$ with $\bar{V}\subseteq\bigcup_{j=1}^{r}U_{j}$ there is $\delta>0$ such that 
   	$$\H_{V}\cap\tv_{\delta}\subseteq\mathfrak{h}^{c}_{U_{1}}\cdots\mathfrak{h}^{c}_{U_{r}}.$$
   	If in the situation above for any $\epsilon>0$ there is $\delta>0$ such that 
   	$$\H_{V}\cap\tv_{\delta}\subseteq(\mathfrak{h}^{c}_{U_{1}}\cap\tv_{\epsilon})\cdots(\mathfrak{h}^{c}_{U_{r}}\cap\tv_{\epsilon}),$$
   	then we say that $\mh$ satisfies the local fragmentation property for $\tau_{co}$. 
   \end{definition}
   
   \begin{lemma}
   	\label{l: weak to strong local fragmentation} For any compact manifold $X$ and $\mh\leq\H(X)$ if $\mh$ satisfies the weak local fragmentation property with respect to $\tau_{co}$ then it satisfies the local fragmentation property with respect to $\tau_{co}$. 
   \end{lemma}
   \begin{proof}
   	Let $V,U_{1},\dots U_{r}\subseteq X$ be open subsets with $\bar{V}\subseteq\bigcup_{i=1}^{r}U_{i}$ and $\epsilon>0$. 
   	Let $N=N(X)$ be as in Fact \ref{f: colorable coverings}. Then $\bar{V}\subseteq\bigcup_{i=1}^{r}\bigcup_{j=1}^{N}V_{i,j}$, where each $V_{i,j}$ is the union of the interiors of some disjoint collection of compact subsets of $V_{i}$ of diameter at most $\frac{\epsilon}{N}$. 
   	The desired conclusion now follows from the weak local fragmentation property applied to $V$ and $\{V_{i,j}\,|\,1\leq i\leq r,1\leq j\leq N\}$ 
   	and the observation that $\mh_{V_{i,1}}\cdots\mh_{V_{i,N}}\subseteq(\tv_{\frac{\epsilon}{N}})^{N}\subseteq\tv_{\epsilon}$.  
   \end{proof}
   
   We recall the following fact, due to Edwards and Kirby \cite{edwards1971deformations}.
   \begin{fact}
   	\label{f: fragmentation} Let $M$ be a compact manifold. Then $\mathcal{H}(X)$ and $\mathcal{H}(X,\partial X)$ have the weak local fragmentation property for $\tau_{co}$.      
   \end{fact}

  \begin{lemma}
  	\label{l: fragmentation}If $(X,\mu)$ is an OU pair, then $\gg:=\hh(X,\mu)$ satisfies the weak local fragmentation property for $\tau_{co}$. 
  \end{lemma}
	  \begin{proof}
	  	Write $\H$ for $\mathcal{H}(X)$ and let $U_{1},\dots U_{r},V$ be open sets in $X$ such that $\overline{V}\subseteq\bigcup_{j=1}^{r}U_{j}$. We will	show by induction on $r$ the existence of some $\delta>0$ such that: 	 
	  	\begin{equation}\label{fragmentation goal}
	  		\gg_{V}\cap\V_{\delta}\subseteq\gg_{U_{1}}^{c}\cdots\gg_{U_{r}}^{c}.
	  	\end{equation}
	  	The weak local fragmentation property for $\tau_{co}$ will then follow. The base case $r=1$ is trivial, so assume $r\geq 2$. 
	    It is easy to find an open set $U'$ such that $\overline{U'}\subseteq\bigcup_{l\neq r}U_{l}$ 
	    and $\overline{V}\subseteq  U'\cup U_{r}$. By the induction hypothesis, there exists $\delta_{0}>0$ such that   
	  	\begin{equation}
	  		 \label{e: fragmentation induction} \gg_{U'}\cap\V_{\delta_{0}}\subseteq \gg_{U_{1}}^{c}\cdots\gg^{c}_{U_{r-1}}.
	  	\end{equation}
	    Next, let $\epsilon:=\frac{\delta_{0}}{2}$. By the local fragmentation property for the ambient group $\H$ there exists $\delta>0$ such that
	  	 \begin{equation}
	  	 	\label{e: fragmentation h}\H_{V}\cap \tv_{\delta}\subseteq (\H_{U'}^{c}\cap\tv_{\epsilon})(\H_{U_{r}}^{c}\cap\tv_{\epsilon}).
	  	 \end{equation}
	  	
	    We claim that condition \eqref{fragmentation goal} holds for the $\delta$ given above. 
	    To see this, choose an arbitrary $g\in\gg_{V}\cap\V_{\delta}$. By \eqref{e: fragmentation h} we can write $g=\hat{g}_{1}\hat{g}_{2}$, where	$\hat{g}_{1}\in\H_{U'}^{c}\cap\tv_{\epsilon}$ and $\quad\hat{g}_{2}\in\H_{U_{r}}^{c}\cap\tv_{\epsilon}$.

	  	By applying \cref{l: ac correction} with $F :=\overline{\supp}(\hat{g}_{1})\cap\overline{\supp}(\hat{g}_{2})$ and $U:=U'\cap U_{r}$, there exists some $h\in\mathcal{H}(X,\partial X)$ such that:
	  	\begin{itemize}
	  		\item $F\subseteq\supp(h)\subseteq\overline{\supp}(h)\subseteq U'\cap U_{r}$,
	  		\item $h\in\tv_{\epsilon}$,
        \item and $\tilde{g}:=\hat{g}_{1}h$ is measure-class preserving on $\supp(h)$.
	  	\end{itemize}
  	  Observe that 
  	  \begin{equation*}
  	  	\overline{\supp}(\tilde{g})\subseteq\supp(\hat{g}_{1})\cup\supp(h)\subseteq\overline{\supp}(\hat{g}_{1})\cup\overline{\supp}(h)\subseteq U'.
  	  \end{equation*}
  	  
  	  Clearly, $\tilde{g}$ agrees with $\hat{g}_{1}$ on $\supp(\tilde{g})\setminus\supp(h)$.
  	  We also have 
  	   \begin{equation}
  	   	\supp(\tilde{g})\setminus\supp(h)\subseteq\supp(\hat{g}_{1})\setminus\supp(h)\subseteq\supp(\hat{g}_{1})\setminus F.
  	   \end{equation}
  	  But $\hat{g}_{1}$ agrees with $g$ on the set $\supp(\hat{g}_{1})\setminus F$. Since $g$ is measure-class preserving everywhere, this implies that $\tilde{g}$ is also measure-class preserving on $X\setminus\supp(h)$. From \cref{o: pasting}, we conclude that $\tilde{g}\in \gg$. The choice of $\epsilon$ also guaranties that $\tilde{g}\in\V_{2\epsilon}\subseteq\V_{\delta_{0}}$. Therefore, by \eqref{e: fragmentation induction} there are $g_{i}\in\gg_{U_{i}}^{c}$, $1\leq i\leq r-1$ such that 
	  	$\tilde{g}=g_{1}\cdots g_{r-1}$.

	  	Finally, write $g_{r}:=h^{-1}\hat{g}_{2}$. Clearly $g_{2}\in\gg$, since $\tilde{g}$ and $g$ are, while we also have 
	  	\begin{equation*}
	  		\overline{\supp}(g_{2})\subseteq \overline{\supp}(\hat{g}_{2})\cup\overline{\supp}(h)\subseteq U_{r}.
	  	\end{equation*}
	  \end{proof}

    The following result follows from an approximation argument such as Lemma 5 in \cite{ihli2022generation}. 
    \begin{corollary}
    	\label{c: fragmentation finer} Let $X=[0,1]$ and suppose that $U_{1},\dots U_{r}$ are open sets with $(0,1)=\bigcup_{j=1}^{r}U_{j}$. Then there exists some continuous function $\epsilon:X\to\R_{\geq 0}$ with $\epsilon^{-1}(0)=\{0,1\}$ such that  
    	$$\{g\in\gg\,|\,\forall x\in V\,\,d(g(x),x)\leq\epsilon(x)\}\subseteq\gg_{U_{1}}^{c}\cdots\gg_{U_{r}}^{c}.$$  
    \end{corollary}

 \subsection*{Uniform perfectness}

   \begin{lemma}
   	\label{l: translation}Assume that $D,D'\subseteq X$ are either $\mu$-adapted balls or boundary half-balls in $X$ satisfying $D\subseteq \mathring{D'}$. In case they are boundary half-undefined()balls, assume in the latter case that $dim(X)\geq 2$. Then there is  
   	$h\in\gg_{D'}$ such that $\{h^{k}(D)\}_{k\in\Z}$ are disjoint. 
   \end{lemma}
   \begin{proof}
   	Assume without loss of generality that $D$ and $D'$ are embedded balls and pick some $\mu$-adapted ball $D_{0}$ with $D\subseteq\mathring{D_{0}}\subseteq D_{0}\subseteq\mathring{D'}$ and some point $p_{0}\in \mathring{D'}\setminus D_{0}$. Using Fact \ref{f: OU uniqueness}, one can easily see that $D_{0}$ is the first entry of a sequence of $\mu$-adapted balls $\{D_{k}\}_{k\geq 0}$ such that for all $k\geq 0$ 
   	\begin{itemize}
   		\item $D_{l+1}\subset\mathring{D_{l}}$ for all $0\leq l<k$,
   		\item there is some constant $\lambda\in(0,1)$ such that $\mu(D_{k+1})=\lambda\mu(D_{k})$ for all $k\geq 0$, 
   		\item 
   		and $D_{k}$ converges to $\{p_{0}\}$ in the Hausdorff topology as $k$ goes to infinity.
   	\end{itemize}
   	Using Observation \ref{o: ball correction for any} it is easy to construct a sequence of elements $(h_{k})_{k\geq 0}\subset\mathcal{H}(X)$ supported on 
   	$D'$ such that
   	\begin{itemize}
   		\item $h_{k}$ agrees with $h_{k+1}$ on $D'\setminus\mathring{D_{k}}$,
   		\item $h_{k}$ maps $D_{l}$ homeomorphically onto $D_{l+1}$ for $0\leq l\leq k$.  
   	\end{itemize}
    In view of the statement of \ref{o: ball correction for any}, we may also assume:
    \begin{align*}
   	\nn{h_{k}}_{bac,D'\setminus D_{k}}=&(\mu(D'\setminus D_{1})-\mu(D'\setminus D_{0}))+\sum_{l=0}^{k-1}(\mu(D_{l})-\mu(D_{l+1})-(\mu(D_{l+1})-\mu(D_{l+2}))\leq\\
   	&\mu(D_{0}\setminus D_{1})+\sum_{l=0}^{k-1}(1-\lambda)^{2}\lambda^{l}\mu(D_{0})\leq\mu(D_{0}\setminus D_{1})+(1-\lambda)\mu(D_{0}).
    \end{align*}
   	The sequence $(h_{k})_{k\geq 0}$ can easily be seen to converge in $\tau_{ac}$ to some $h_{\infty}\in\gg$ with the desired property.  	
   \end{proof}
  
   The following follows easily from a standard argument due to Anderson \cite{anderson1958algebraic} and is likely known.
  \begin{corollary}
  	\label{c: uniform perfectness} Let $D,D'\subseteq X$ be two embedded (boundary) balls such that $D\subseteq \mathring{D'}$. Then for any 
  	$g\in\gg_{\mathring{D}}$ there exists $h_{1},h_{2}\in\gg_{\mathring{D'}}$ such that $g=[h_{1},h_{2}]$. 
  \end{corollary} 
  \begin{proof}
   Take $h_{2}$ as given by Lemma \ref{l: translation} and let $h_{1}\in\gg_{\bigcup_{k\geq 0}h^{k}(\mathring{D})}$ be the element that acts on $h^{k}(\mathring{D})$ as $g^{h^{-k}}$ for all $k\geq 0$ (using Observation \ref{o: pasting}). 
  \end{proof}

\section{Some general preliminaries for Theorem \ref{t: main}}
  \label{s:preliminaries}
  
  Here we collect some general definitions and facts that will be needed later. Most of it applies in the general setting of Proposition \ref{p: generalization}.

 \subsection*{Martingales and Chebyshev's inequality}

\begin{definition}
	\label{l: weakly independent} Let $(X,\mathcal{A},\mu)$ be a probability space. We say that a sequence of functions $\psi_{1},\dots \psi_{n}:X\to\R$ in $L^{1}(X,\mu)$ is a sequence of martingale increments 
	if as random variables the expectation of $\psi_{k}$ conditional on $(\psi_{l})_{l<k}$ is $0$. In other words, if $\int\psi_{1}\dm=0$ and for every $1\leq k\leq n-1$ and every finite tuple of real intervals
	$\underline{J}:=(J_{1},\dots J_{k})$ if we let $A_{\underline{J}}:=\{x\,|\,\psi_{i}(x)\in J_{i}\,\text{for $1\leq i\leq k$}\}$ then we have 	$\int_{A_{\underline{J}}}\psi_{k+1}\dm=0$. \footnote{Equivalently, if the sequence $(\sum_{j=1}^{k}\psi_{j})_{k=1}^{n}$ is a martingale that starts at zero.}
\end{definition}

For the following, see, for instance, \cite{klenke2013probability}, Example 10.2 and Theorem 10.4. 
\begin{fact}
	\label{f: square variation}For any sequence of martingale increments $\psi_{1},\dots \psi_{n}$ we have $\nn{\sum_{l=1}^{n}\psi_{l}}^{2}_{2}=\sum_{l=0}^{n}\nn{\psi_{l}}^{2}_{2}$. 
\end{fact}

\begin{fact}[Chebyshev's inequality]
  	\label{f: Chebyshev's inequality} For any probability space $(X,\mathcal{A},\mu)$, any $f\in L^{2}(X,\mu)$ and any $\epsilon>0$ we have 
  	  \begin{equation*}
  	  	\mu(\{x\in X\,|\,|f(x)|>\epsilon\})\leq\frac{\nn{f}_{2}^{2}}{\epsilon^{2}}.
  	  \end{equation*}	
\end{fact}

The following is a standard fact and an immediate corollary of Facts \ref{f: Chebyshev's inequality} and \ref{f: square variation}.  
\begin{corollary}
	\label{c: Chebyshev martingale}For any $\epsilon,C>0$ there exists some $N=N(C,\epsilon)\in\N$ such that for any $n\geq N$	and any sequence of martingale increments $\psi_{0},\psi_{1},\dots$ with $\nn{\psi_{l}}^{2}_{2}\leq C$ for all $1\leq l\leq n$ we have 
	$$
	  \mu(\{\,x\,|\,|\frac{\sum_{l=0}^{n}\psi_{l}(x)}{n}|>\epsilon\})<\epsilon.
	$$
\end{corollary}

  \subsection*{Two observations on derivatives and integrals} 

 Let $X$ and $\mu$ be as in the assumptions of Proposition \ref{p: generalization} and write $\gg:=\hh(X,\mu)$. Without loss of generality, let $\mu$ be a probability measure. Given a real function $f$ on $X$, we write $\phi_{+}:=\phi\cf_{\{\phi\geq 0\}}$.

\begin{lemma}
	\label{l: l1bound}Let $g\in\gg$ and let $A\in \mathcal{B}(X)$. Then 
	$$\nn{g}_{ac,A}=2\int_{A}(1-\mathcal{D}(g))_{+}\dm+(\ms{g(A)}-\ms{A})\leq \ms{A}+\ms{g(A)}.$$
\end{lemma} 
\begin{proof}
	To begin with, note that
	$$\mu(g(A))-\mu(A)=\int_{A}(\mathcal{D}(g)-1)\dm=\int_{A}(\mathcal{D}(g)-1)_{+}\dm-\int_{A}(1-\mathcal{D}(g))_{+}\dm.$$
	Since $f_{+}+(-f)_{+}=|f|$ for all $f\in L^{1}(X,\mu)$, it follows that
	\begin{equation*}
		2\int_{A}(1-\mathcal{D}(g))_{+}\dm+\mu(g(A))-\mu(A)=\int_{A}|1-\mathcal{D}(g)|\dm=\nn{g}_{ac,A}.
	\end{equation*}
	On the other hand, since $\mathcal{D}(g)>0$ almost everywhere, we have 
	$\int_{A}(1-\mathcal{D}(g))_{+}\dm\leq \ms{A}$, from which the inequality in the statement follows. 
\end{proof}

\begin{corollary}
	\label{c: often enough negative enough} Suppose that $A\in\mathcal{B}(X$) and $g\in\gg$ satisfy
	\begin{enumerate}
	\item \label{cond1}$\nn{g}_{ac,A}\geq\epsilon\ms{A}$ for some $\epsilon>0$,
	\item \label{cond2} and $\ms{g(A)}\leq(1+\delta)\ms{A}$ for some $\delta>0$.
	\end{enumerate}
	 Then 
	 $$\mu(A_{0})\geq\frac{(\epsilon-\delta)\ms{A}}{4},$$ 
	 where $A_{0}:=\{x\in A\,|\,\mathcal{D}(g)(x)\leq 1-\frac{\epsilon-\delta}{4}\}$.
\end{corollary}
\begin{proof}
	By Lemma \ref{l: l1bound} we have 
	$$
	 \int_{A}(1-\mathcal{D}(g))_{+}\dm=\frac{\nn{g}_{ac,A}-(\ms{g(A)}-\ms{A})}{2}\geq \frac{\epsilon-\delta}{2}\ms{A}.
	$$
  On the other hand, we have:
	$$
	\int_{A}(1-\mathcal{D}(g))_{+}\dm\leq\mu(A_{0})+\frac{\epsilon-\delta}{4}\mu(A\setminus A_{0})\leq\mu(A_{0})+\frac{\epsilon-\delta}{4}\mu(A),
	$$
	since $\mathcal{D}(g)> 0$ a.e.. The result follows.
\end{proof}

 \subsection*{Filling blocks and good covering sequences}
 
 \begin{definition}
 	\label{d: filling set} Given $\epsilon>0$ we say that $A\subseteq X$ is an $\epsilon$-filling block if $\mu(Fr(A))=0$, $\mu(A)>0$ and there exist $n>1$ and  a finite set $\{g_{l}\}_{l=1}^{n}\subseteq\gg$ such that $\mathcal{D}(g_{l})$ is constant almost everywhere in some neighborhood of $\bar{A}$ and 
  if we let $A_{l}:=g_{l}(A)$ we have: 
  \begin{itemize}
	  	\item $\mu(A_{l}\cap A_{l'})=0$ for $1\leq l<l'\leq n$
	  	\item $\mu(X\setminus\bigcup_{l=1}^{n}A_{l})<\epsilon$. 
  \end{itemize} 
   \end{definition}
 
 The following result is due to Fathi~\cite{fathi1980structure}, based on the work of Brown for the unmeasured case; see also \cite[Chapter 17]{CV1977}.  
 
 \begin{fact}\label{lem:brown-measure}
 	For each OU pair $(X,\mu)$ there exists a continuous surjection \[f\colon \bar{B}^n(1)
 	\to X\] such that the following hold:
 	\begin{enumerate}[label=(\alph*), ref=(\roman{enumi}\alph*)]
 		\item 	the restriction of $f$ to $B^{n}(1)=\Int(\bar{B}^{n}(1))$ is an embedding onto an open dense subset of $X$, so that in particular, $f(B^{n}(1))\subseteq\Int(X)$ and thus $\partial X\subseteq  f(\partial \bar{B}^n(1))$;
 		\item 	we have $f(B^{n}(1)))\cap f(\partial \bar{B}^n(1))=\varnothing$; 
 		\item\label{p:brown-leb} the measure $\mu$ is the push-forward of the Lebesgue measure on $\bar{B}^{n}(1)$ by $f$.
 	\end{enumerate}
 \end{fact}

    \newcommand{\lk}[0]{\mathscr{L}}
   \newcommand{\per}[0]{\mathfrak{p}}
  \begin{definition}
   	   \label{d: partitions}By a \emph{covering sequence} (for $(X,\mu)$) we mean a sequence $$\mathscr{S}=(\mathcal{P}_{m},\lk_{m},\lk'_{m})_{m\geq 0}$$ where $\mathcal{P}_{m}=\{P(m,l)\}_{l\in\lk_{m}}$ is a collection of  finitely many measurable sets and $\lk'_{m}\subseteq\lk_{m}$ is such that the following is satisfied:
   	   \begin{enumerate}[label=(\arabic*)]
   	   	\item \label{en-2} $\bigcup_{l\in\lk_{m}}P(m,l)=X$, 
   	   	\item \label{en-1} for distinct $l,l'\in\lk_{m}$ we have $P(m,l)\cap P(m,l')\subseteq\Fr(P(m,l))\cap\Fr(P(m,l'))$,  
   	   	\item \label{en0} $\mu(\Fr(P(m,l)))=0$ and $\mu(P(m,l))>0$ for all $m\geq 0$ and $l\in\lk_{m}$,
   	   	\item \label{en1}	for any $m\geq 0$ and $l\in\lk_{m+1}$ there is $l'\in\lk_{m}$ with $P(m,l)\subseteq P(m+1,l')$,
 	      \item \label{en2} $\sup_{l\in\lk_{m}}\diam(P(m,l))$ converges to $0$ as $m$ tends to infinity,
   	   	\item \label{en4} for $m'>m\geq 0$ and every $l\in\lk_{m}$ we have $$\frac{1}{2}\mu(P(m,l))\leq \mu(P(m,l)\cap\bigcup_{l\in\lk'_{m'}}P(m',l')).$$  
   	   \end{enumerate}   	     
       If in addition to this for some $\epsilon$ we have 
       \begin{enumerate}[label=(\arabic*)]  
       	\setcounter{enumi}{5}
       	\item \label{en5} $P(m,l)$ is an $\epsilon$-filling block for all $m\geq 0$, $l\in\lk'_{m}$,
       \end{enumerate}
       then we say $\mathcal{P}$ is \emph{$\epsilon$-filling}.
       
       We say that $(X,\mu)$ is \emph{well-covered} if it admits an $\epsilon$-filling good covering sequence for all $\epsilon>0$. 
   \end{definition}

 \begin{corollary}
 	\label{c: ou filling} If $(X,\mu)$ is an OU pair, then any $\mu$-adapted ball $D\subseteq X$ is an $\epsilon$-filling block
 	for all $\epsilon>0$ and $(X,\mu)$ is well-covered.
 \end{corollary}
 \begin{proof}
 	That any $\mu$-adapted ball is an $\epsilon$-filling block follows from Fact \ref{lem:brown-measure} and Corollary \ref{o: ball transitivity}. To construct $\mathscr{S}$ as above, take the continuous surjection $f:I^{n}\cong B_{1}\to X$ given by \ref{lem:brown-measure} and for each $m\geq 0$ push onto $X$ the cover of $I^{n}$ by the collection of all cubes of the form $[0,\frac{1}{2^{m}}]^{m}+\bar{\alpha}$, where $\bar{\alpha}$ is an arbitrary $n$-tuple with values in $\frac{1}{2^{m}}(\N\cap[0,2^{m}-1])$.
 \end{proof}

\section{Density of elements with well-distributed norm}
\label{s:elements spread}

\label{s: well distributed elements}
   Throughout this section and the next we assume that $(X,\mu)$ satisfies the hypotheses of Proposition \ref{p: generalization}, with $\mu(X)=1$, and $d$ is a compatible metric on $X$.

   \newcommand{\qer}[0]{\mathfrak{q}}
   \begin{definition}
	   \label{d: perturbation} Given $A\in\mathcal{B}(X)$ with $\mu(Fr(A))=0$, $\mu(A)>0$ and $\epsilon>0$ we write
	    \begin{equation*}
	    	\per_{A}(\epsilon)=\sup \{\delta>0\,|\,\mu(g(A)\triangle A)<\epsilon\mu(A)\,\},
	    \end{equation*}
	    so that, in particular, $|\mu(A)-\mu(g(A))|<\epsilon$ if $g\in\V_{\delta}$ for $\delta<\per_{A}(\epsilon)$. 
   \end{definition}
   For example, if $X=I$, $\mu$ the Lebesgue measure, $d$ the standard distance, and $A\subseteq\mathring{I}$ an interval, then $\per_{A}(\epsilon)=\frac{\epsilon}{2\mu(A)}$ for $\epsilon$ small enough. 
   
   \begin{observation}
   	\label{o: perturbation}$\per_{A}(\epsilon)>0$ for $\epsilon>0$ and $A$ as in \ref{d: perturbation}. 
   \end{observation}
   \begin{proof}
   	It suffices to show that for any sequence $(g_{n})_{n}\in\mathcal{H}(X)$ that converges to $\Id$ in $\tau_{co}$ the measure of $\mu(g_{n}(A)\setminus A)$ and that of $\mu(A\setminus g_{n}(A))$ both converge to $\mu(A)$. As in the proof of Lemma \ref{l: completeness}, 
   	we must have $\mathring{A}\subseteq\bigcup_{n\in\N}\bigcap_{m\geq n}g_{n}(A)$ and thus $\bigcap_{n\in\N}\bigcup_{m\geq n}(A\setminus g_{n}(A))\subseteq\Fr(A)$, which implies that 
   	$\lim_{n}\mu(g_{n}(A)\setminus A)=0$. Arguing with $X\setminus A$ instead of $A$ one can use the same argument to deduce that $\lim_{n}\mu(g_{n}(A)\setminus A)=0$.
   \end{proof}
   
   We call to the reader's attention that the proof of the following Lemma contains a tree of sublemmas. As we will do in the proof of \cref{p: almost independent}, we use the shading of the square closing the proof to indicate their degree of nesting.
   \begin{lemma}
   	\label{l: independent increments} Given any $g_{0}\in\gg$, any compact $A\subseteq X$ and any $\eta>0$ there is some $\delta(g_{0},\eta)>0$ 
   	such that for any $g_{1}\in\V_{\delta}$ we have 
   	$$\nn{g_{0}\circ g_{1}}_{ac,A}\geq\max\{\nn{g_{0}}_{ac,A},\frac{1}{8}\nn{g_{1}}_{ac,A}\}-\eta.$$ 
   \end{lemma}
   \begin{proof}
   	We start by choosing $\phi\in C^{0}(X,\mathbb{R}_{> 0})$ such that $\nm{\mathcal{D}(g_{0})-\phi}\leq\frac{\eta}{4}$ 
    and $\lambda>0$ such that $\omega_{\phi}(2\lambda)<\frac{\eta}{8}$.
   Consider now some partition $A=\coprod_{i=1}^{n}A_{i}$ into finitely many Borel sets of diameter less than $\lambda$ and with a $\mu$-null frontier.
   	 For $1\leq i\leq n$ choose an arbitrary $x_{i}\in A_{i}$.  Let $\mathcal{I}_{+}$ be the set of indices $i$ such that $\phi(x_{i})>1$ and $\mathcal{I}_{-}:=\{1,\dots n\}\setminus\mathcal{I}_{+}$. 
   	 We will now check that any 
   	 $\delta$ with
   	 \begin{equation*}
   	 	0<\delta<\min \,\{\,\lambda,\,\per_{A}(1),\per_{A_{i}}(2),\,\per_{A_{i}}(\frac{\eta}{4n\phi(x_{i})\mu(A_{i})})\,\}_{i=1}^{n}
   	 \end{equation*}
   	 verifies the conditions of the statement.    	 
   	 Fix some arbitrary $g_{1}\in\V_{\delta}$ and let $\psi:=\sum_{i=1}^{n}\phi(x_{i})\ind_{A_{i}}$.
   \begin{lemma}
   	\label{l: first inequality}The following inequality holds: 
   	\begin{equation}
   		  \label{eq:first inequality}\nn{\psi \mathcal{D}(g_{1})-1}_{1,A}\geq \max\,\{\, \nn{\psi-1}_{1,A},\,\frac{1}{8}\nn{g_{1}}_{ac,A}\,\}-\frac{\eta}{4}.
   	\end{equation}
   \end{lemma}	
 \begin{subproof}
 	    We will show that for any $i\in\{1,2,\dots n\}$ we have:
   	  \begin{equation}
   	  	\label{e: interval equation} \int_{A_{i}}|\phi(x_{i})\mathcal{D}(g_{1})-1|\dm\geq \max\,\{\, |\phi(x_{i})-1|\mu(A_{i}),\frac{1}{8}\nn{g_{1}}_{ac,A_{i}}\, \}-\frac{\eta}{4n}.
   	  \end{equation}
   	  The inequality in the statement can be obtained by adding the one above for all values $1\leq i\leq n$.    For any $i\in\mathcal{I}^{+}$ using the fact that $g_{1}\in\V_{\delta}$ and $\delta<\frac{\eta}{4n\phi(x_{i})\mu(A_{i})}$ we get: 
      \begin{align}
      \label{basic inequality}
      \begin{split}
	      \int_{A_{i}}|\phi(x_{i})\mathcal{D}(g_{1})-1|\dm\geq& \int_{A_{i}}(\phi(x_{i})\mathcal{D}(g_{1})-1)\dm =\phi(x_{i})\int_{A_{i}}\mathcal{D}(g_{1})\dm-\ms{A_{i}}=\\
  	   	\phi(x_{i})\mu(g_{1}(A_{i}))-\ms{A_{i}}\geq& \phi(x_{i})(\ms{A_{i}}-\frac{\eta}{4n\phi(x_{i})})-\ms{A_{i}}=|\phi(x_{i})-1|\mu(A_{i})-\frac{\eta}{4n}.
  	  \end{split}
     \end{align}
    For $i\in\mathcal{I}_{-}$ we can argue similarly:
     \begin{align}
     \label{basic inequality 2}
      \begin{split}
	      &\int_{A_{i}}|\phi(x_{i})\mathcal{D}(g_{1})-1|\dm\geq  \int_{A_{i}}(1-\phi(x_{i})\mathcal{D}(g_{1}))\dm =\ms{A_{i}}-\phi(x_{i})\mu(g_{1}(A_{i}))\geq\\
  	   	& \ms{A_{i}}-\phi(x_{i})(\ms{A_{i}}+\frac{\eta}{4n\phi(x_{i})})=(1-\phi(x_{i}))\mu(A_{i})-\frac{\eta}{4n}=|\phi(x_{i})-1|\mu(A_{i})-\frac{\eta}{4n}.
  	  \end{split}
     \end{align}
      
    This establishes the first half of \eqref{e: interval equation} for an arbitrary value $i\in\{1,\dots n\}$. It remains to show 
    \begin{equation*}
    	\int_{A_{i}}|\phi(x_{i})\mathcal{D}(g_{1})-1|\dm\geq\frac{1}{8}\nn{g_{1}}_{ac,A_{i}}-\frac{\eta}{4n},
    \end{equation*}
    and, clearly, we may assume
    \begin{equation}
    	\label{e: assumption}|\phi(x_{i})-1|\ms{A_{i}}<\frac{1}{8}\nn{g_{1}}_{ac,A_{i}},
    \end{equation}
     so that then  
    \begin{equation*}
    	\label{e: small g0}|\phi(x_{i})-1|<\frac{1}{8\ms{A_{i}}}\nn{g_{1}}_{ac,A_{i}}\leq\frac{1}{8\ms{A_{i}}}(\ms{A_{i}}+\ms{g_{1}(A_{i})})\leq\frac{4\ms{A_{i}})}{8\ms{A_{i}}}\leq\frac{1}{2},
    \end{equation*}
    where we have used Lemma \ref{l: l1bound} in the second inequality and the fact that $g_{1}\in\V_{\delta}$ and $\delta<\per_{A_{i}}(2)$ in the last one. It follows that $\phi(x_{i})\geq\frac{1}{2}$.
      
      We will also need the following:    
   		\begin{lemma}
   		\label{l: bound}Let $A\in\mathcal{B}(X)$ and $\theta\in L^{1}(X,\mu)$. 
   		Then for all $\beta>0$ we have:  		 
   		 \begin{equation*}
   		 	\label{basiineq}\nn{\beta \theta-1}_{1,A}\geq \beta\nn{\theta-1}_{1,A}-|\beta-1|\ms{A}.
   		 \end{equation*}
   	\end{lemma}
   	\begin{subsubproof}
   		Let $A_{+}:=\{x\in A\,|\,\theta(x)\geq 1\}$ and $A_{-}:=\{x\in A\,|\,\theta(x)<1\}$.    
   	  We also have the inequality
   		\begin{align*}
   		\begin{split}
   			\int_{A_{-}}|\beta \theta-1|\dm\geq
   		\int_{A_{-}}\beta(1-\theta)\dm+\int_{A_{-}}(1-\beta)\dm&=\beta\nn{\theta-1}_{1,A_{-}}\dm+(1-\beta)\mu(A_{-}).
   		\end{split}
   		\end{align*}
   		Arguing analogously, one can also obtain:  
   		\begin{align*}
   			\int_{A_{+}}|\beta \theta-1|\dm\geq \int_{A_{+}}\beta(\theta-1)\dm+\int_{A_{+}}(\beta-1)\dm=\beta\nn{\theta-1}_{ac,A_{+}}+(\beta-1)\mu(A_{+}).
   		\end{align*}
   		Adding the two inequalities above yields the desired result, since clearly $$(1-\beta)(\mu(A_{-})-\mu(A_{+}))\leq|1-\beta|\mu(A).$$   	
   	\end{subsubproof}

   	
   	Applying Lemma \ref{l: bound} with $\beta=\phi(x_{i})$ for the first inequality below, and then combining the inequalities $\phi(x_{i})\geq\frac{1}{2}$ and  \eqref{e: assumption} for the second one, we get: 
    $$
   	\nn{\phi(x_{i})\mathcal{D}(g_{1})-1}_{ac,A_{i}}\geq \phi(x_{i})\nn{g_{1}}_{ac,A_{i}}-|\phi(x_{i})-1|\ms{A_{i}}\geq\frac{1}{2}\nn{g_{1}}_{ac,A_{i}}-\frac{1}{8}\nn{g_{1}}_{ac,A_{i}}\geq\frac{1}{8}\nn{g_{1}}_{ac,A_{i}},
   	$$
   	ending the proof of \cref{l: first inequality}. 
   	\end{subproof}

   	To conclude the proof of Lemma \ref{l: independent increments} we apply \cref{o: approximation use}. We use the fact that $\nn{\mathcal{D}(g_{0})-\phi}_{1,A}\leq\frac{\eta}{4}$, as well as $\delta<\lambda$, $g_{1}\in\V_{\delta}$, and $\omega_{\phi}(2\lambda)<\frac{\eta}{8}$. We obtain:
      \begin{equation}
    	\label{e: first main ineq} \nn{g_{0}\circ g_{1}}_{ac,A}\geq \nn{\psi \mathcal{D}(g_{1})-1}_{1,A}-\frac{\eta}{4}-\frac{\mu(g_{1}(A))\eta}{8}\geq \nn{\psi \mathcal{D}(g_{1})-1}_{1,A}-\frac{3\eta}{8}.
    \end{equation}
    On the other hand, notice that by the choice of $\lambda$ and $\phi$ we have 
    $$\nn{\psi-1}_{1,A}\geq\nn{\phi-1}_{1,A}-\frac{\eta}{8}\geq\nn{g_{0}}_{ac,A}-\frac{3\eta}{8},$$
    which together with Lemma \ref{l: first inequality} yields
    \begin{align*}
    	\nn{\psi \mathcal{D}(g_{1})}_{1,A}-\frac{3\eta}{8}\geq \max\,\{\nn{\psi-1}_{1,A},\frac{1}{8}\nn{\mathcal{D}(g_{1})}_{ac,A}\}-\frac{5\eta}{8}\geq \max\,\{\nn{g_{0}}_{ac,A},\frac{1}{8}\nn{\mathcal{D}(g_{1})}_{ac,A}\}-\eta.
    \end{align*}
    Combining this with \eqref{e: first main ineq}, the conclusion of \cref{l: independent increments} follows.
   \end{proof} 

   \begin{lemma}
	  \label{l: derivative everywhere} Suppose that $\tau$ is a group topology on $\gg$ finer than $\tau_{co}$ and such that there is $\epsilon_{0}>0$ with the property that for all $\V\in\idn$ there is $g\in\V$ with
	  $\nn{g}_{ac}>\epsilon_{0}$. 
	  
	  Then for every $\V\in\idn$ any $\frac{\epsilon_{0}}{4}$-filling block $A$ there is some $g\in\V$ such that $\nn{g}_{ac,A}>\frac{\epsilon_{0}}{4}\mu(A)$. 
   \end{lemma}
   \begin{proof}
   	 Assume for the sake of contradiction the existence of some $\V\in\idn$ and an $\frac{\epsilon_{0}}{4}$-filling block $A$, such that for all $g\in\mathcal{V}$ we have $\nn{g}_{ac,A}\leq\frac{\epsilon_{0}}{4}\mu(A)$.  
   	 
   	 According to \cref{d: filling set} there exists $n>1$ and a finite set $\{g_{l}\}_{l=1}^{n}\subseteq\gg$ such that $\mathcal{D}(g_{l})$ is constant on some neighborhood $U_{l}$ of $A$ and if we let $A_{l}:=g_{l}(A)$ we have: 
     \begin{itemize}
	  	\item $\mu(A_{l}\cap A_{l'})=0$ for $1\leq l<l'\leq n$
	  	\item $\mu(X\setminus\bigcup_{l=1}^{n}A_{l})<\epsilon_{0}$. 
     \end{itemize} 
   	 
   	 If we let $\delta_{l}>0$ be small enough that $\Nb_{\delta_{l}}(A_{l})\subseteq g_{l}(U)$, then for any element $g\in\V_{\delta_{l}}$ a simple application of the chain rules implies that for almost all $x\in A$ we have 
   	 $\mathcal{D}(g^{g_{l}})(x)=\mathcal{D}(g)(g_{l}(x))$, from which the equality 
   	 $\mu(A_{l})\nn{g}_{ac,A}=\mu(A)\nn{g}_{ac,A_{l}}$ follows.
   	 
   	 Let $\delta^{*}=\min_{1\leq l\leq n}\delta_{l}$ and consider $\mathcal{V}':=\V_{\delta^{*}}\cap\bigcap_{l=1}^{n}\V^{g_{l}^{-1}}\in\idn$ (where we use $\tau_{co}\subseteq\tau$). By assumption, we can find $g\in\V'$ with the property that $\nn{g}_{ac}>\epsilon_{0}$. In particular $\nn{g^{g_{l}}}_{ac,A}\leq\frac{\epsilon_{0}}{4}\mu(A)$ for all $1\leq l\leq n$. By the discussion of the previous paragraph  $\nn{g}_{ac,A_{l}}\leq\frac{\epsilon_{0}}{4}\mu(A_{l})$ for all $1\leq l\leq n$, and if we write $C:=\bigcup_{l=1}^{n}A_{l}$, then clearly also  
   	 $\nn{g}_{ac,C}=\frac{\mu(C)}{\mu(A)}\nn{g}_{ac,A}\leq\frac{\epsilon_{0}}{4}\mu(C)\leq\frac{\epsilon_{0}}{4}$. Therefore 
   	 \begin{equation*}
   	 	\mu(g(X\setminus C))-\mu(X\setminus C)=\mu(C)-\mu(g(C))\leq\frac{\epsilon_{0}}{4},
   	 \end{equation*}
   	 and thus $\nn{g}_{ac,X\setminus C}\leq 2\mu(X\setminus C)+\frac{\epsilon_{0}}{4}\leq \frac{3\epsilon_{0}}{4}$ by Lemma \ref{l: l1bound}. This leads to a contradiction:
   	  \begin{equation*}
   	  	\nn{g}_{ac}=\nn{g}_{ac,C}+\nn{g}_{ac,X\setminus C}\leq \frac{\epsilon_{0}}{4}+\frac{3\epsilon_{0}}{4}=\epsilon_{0}. 
   	  \end{equation*}
   \end{proof}
  
  \begin{corollary}
  	\label{c: uniform derivative} Let $\epsilon_{0}>0$ be such that for all $\V\in\idn$ there is $g\in\V$ with $\nn{g}_{ac}>\epsilon_{0}$
  	and let $A_{1},\dots A_{n}$ be a collection of $\frac{\epsilon_{0}}{4}$-filling blocks in $X$ whose pairwise intersection is $\mu$-null. Then for any $\V\in\idn$ there is $g\in\V$ such that 
  	$\nn{g}_{ac,A_{i}}\geq\frac{\epsilon_{0}\mu(A_{i})}{64}$ for all $1\leq i\leq n$.
  \end{corollary}
  \begin{proof}
  	Choose $\mathcal{U}\in\idn$ such that $\mathcal{U}^{n}\subseteq\V$ and iteratively $g_{1},\dots g_{n}\in\gg$ so that writing	$\bar{g}_{0}=\Id$ and $\bar{g}_{i}=g_{1}g_{2}\dots g_{i}$, then given $g_{1},\dots g_{i-1}$, we choose $g_{i}$ as follows: 
  	\begin{enumerate}
  		\item Let $\V_{1}$ stand for $\gg$ if $i=1$, and otherwise for $\V_{\delta}$, where $\delta=\delta(\bar{g}_{i-1},\frac{\epsilon_{0}}{2^{n+7}})$ is given by Lemma \ref{l: independent increments}. 
  		\item Using the fact that $\tc\subseteq\tau$ and Lemma \ref{l: derivative everywhere} one can pick $g_{i}\in\mathcal{U}\cap\V_{1}$ with 
  		$\nn{g_{i}}_{ac,A_{i}}\geq\frac{\epsilon_{0}\mu(A_{i})}{4}$ 
  	\end{enumerate} 
   An iterated application of Lemma \ref{l: independent increments} easily shows that $g:=\bar{g}_{n}\in\V$ satisfies the required properties. 
  \end{proof}

   \section{Covering sequences and the density of elements with very large norm}
   \label{s:dense large elements}
   
   \begin{definition}
   	Let $\Lambda: (0,\infty)\to \R$ be the continuous function defined by: 
   	\begin{equation*}
   		\Lambda(t):=  
   		\begin{cases*}
   		\ln(t) & if $t\geq 1$ \\
   		 t-1        & if $t< 1$
   	  \end{cases*}
   	\end{equation*}
   \end{definition}
  
  \begin{observation}
  	\label{o: small variance after perturbation}Let $\phi:X\to\R_{>0}$ and $\rho:X\to [-1,1]$ be measurable functions, $A\in\mathcal{B}(X)$. Then 
  	$$\nn{\Lambda(\phi)+\rho}_{2,A}^{2}\leq  2\nn{\phi}_{1,A}+6\ms{A}$$.
  \end{observation}
  \begin{proof}
  	Consider an arbitrary value $\gamma\in[-1,1]$. For $y\geq 1$ and we have: 
  	\begin{equation*}
  		 (\Lambda(y)+\gamma)^{2}=(\ln(y)+\gamma)^{2}\leq (\ln(y)+1)^{2}\leq 2(1+\ln(y)+\frac{\ln(y)^{2}}{2})\leq 2e^{\ln(y)}=2y,
  	\end{equation*}
  	while for any $0\leq y\leq 1$ we have 
  	\begin{equation*}
  		(\Lambda(y)+\gamma)^{2}=(y+\gamma-1)^{2}\leq 4.
  	\end{equation*}
    One can use this to conclude:
  	$$
  	\nn{\Lambda(\phi)+\rho}_{2,A}^{2}\leq 2\nn{(\phi-1)_{+}+1}_{1,A}+4\mu(\{x\in X\,|\,\phi(x)\leq 1\})\leq 2\nn{\phi}_{1}+4\ms{A}.
  	$$  	 
  \end{proof}

\newcommand{\nf}[0]{\mathbf{m}}
\newcommand{\df}[0]{\boldsymbol{\delta}}
  \begin{proposition}
  	\label{p: almost independent}For $\epsilon >0$ small enough the following holds. Let $(\mathcal{P}_{m},\lk_{m},\lk'_{m})_{m\geq 0}$ be a good covering sequence as defined in Definition \ref{d: partitions}.
  	For any $\eta>0,k\in\mathbb{N}^{+}$ there are functions 
  	$$ \df_{\eta,k}:\gg^{k-1}\to\R_{>0},\quad\quad\nf_{\eta,k}:\gg^{k}\to\N $$  	
  	with the properties described below. \footnote{Observe that, in particular, $\nf_{1,\eta},\df_{1,\eta}$ are constants dependent on $\eta$.} 
  	
  	For any sequence
  	$g_{1},\dots g_{n}$ of elements of $\gg$, writing $m_{0}=0$ and also   
  	$$\delta_{k}:=\df_{\eta,k-1}(g_{1},\dots ,g_{k-1}), \quad\quad m_{k}:=\nf_{\eta,k}(g_{1},\dots g_{k})$$ 
  	for $k\geq 1$, if the conditions  
   	\begin{enumerate}[label=(\Roman*)]
  		\item \label{assumption 1}$g_{k}\in\V_{\hat{\delta}_{k}}$
  		\item \label{assumption 2}$\nn{g_{k}}_{ac,P(m_{k-1},l)}>\epsilon\mu(P(m_{k-1},l))$ for all $l\in\lk_{m_{k}}'$ 
  	\end{enumerate}
  	hold, then there are 
  	\begin{itemize}
  		\item $(\psi_{k})_{k=1}^{n}$ is some sequence of martingale increments on $(X,\mu)$ such that $\nn{\psi_{k}}_{2}^{2}\leq 8$,
  		\item $\epsilon':=\frac{\epsilon}{64}(-\frac{\epsilon}{32}-\ln(1-\frac{\epsilon}{32}))>0$,
  		\item and $\Theta\in L^{1}(X,\mu)$ with $\nn{\Theta}_{1}\leq\eta$,
  	\end{itemize}
    such that the inequality
  	$$
  	(g_{1}\circ g_{2}\cdots g_{n})'\leq \Theta+ e^{-n\epsilon'+\sum_{k=1}^{n}\psi_{l}}
  	$$ 
  	holds almost everywhere. 
  \end{proposition}
  
  \begin{proof}
  
    Fix $\eta>0$. 
    
    \subsection*{Outline of the proof} The construction is by induction on $k$. 
    We assume that $\df_{j,\eta},\nf_{j,\eta}$ are already defined for $1\leq j<k$ and we are given a tuple $\bar{g}_{k-1}:=(g_{1},\dots g_{k-1})\in\gg^{k-1}$, which satisfies the conditions of the statement of \cref{p: almost independent}. We also assume that the construction has provided us with step functions $\psi_{1},\dots \psi_{k-1}$, together with some additional auxiliary data, to be described in detail below. 
    
    We begin the $k$-th step of the construction by using the previous information to find $\delta_{k}:=\df_{\eta,k}(\bar{g}_{k-1})$.
    After that, given $g_{k}\in\gg$ satisfying
    \ref{assumption 1} and \ref{assumption 2}, we describe how to find the auxiliary data aforementioned, consisting of
    \begin{itemize}
    	\item functions $\theta,\phi:X\to \R_{>0}$ (see below for their properties and the construction), 
    	\item and a positive constant $\lambda_{k}>0$,
    \end{itemize}
     after which we choose $m_{k}:=\nf_{\eta,k}(\bar{g}_{k-1},g_{k})$ and finally construct $\psi_{k}$ using all the previous choices. We conclude by showing that for any tuple $(g_{1},\dots g_{n})\in\gg^{n}$ as in the hypothesis and the resulting sequence of step functions $\psi_{1},\dots\psi_{n}$ obtained in the above process, there is $\Theta\in L^{1}(X,\mu)$ verifying the conditions of the statement.
    
    \subsection*{Notation} We begin by establishing some notation (recall that $m_{0}=0$):
    \begin{itemize}
    	\item $\hat{\epsilon}:=(-\frac{\epsilon}{32}-\ln(1-\frac{\epsilon}{32}))$ so that $\epsilon'=\frac{\epsilon\hat{\epsilon}}{64}$,
    	\item $\Delta_{j}:=e^{-j\epsilon'+\sum_{l=1}^{j}\psi_{l}}$ (in particular, $\Delta_{0}=1$),
    	\item $\per_{j}:=\min_{l\in\lk_{m_{j}}}\per_{P(m_{j},l)}$ (see Definition \ref{d: perturbation}),
    	\item $\nu_{j,l}:=\mu(P(m_{j},l))$.
    \end{itemize}
    
    \subsection*{Choice of $\delta_{k},m_{k}$ and the auxiliary data} To begin with, take (adopting the convention $\lambda_{0}:=1$):
    \begin{equation}
    	\label{e: delta choice} \delta_{k}:=\frac{1}{2}\min \{\,\lambda_{k-1},\per_{k-1}(\frac{\epsilon'}{2}),\per_{k-1}(\frac{\epsilon}{4}),\delta_{k-1}\,\})
    \end{equation} 
     
    Given $g_{k}\in\gg$ satisfying conditions \ref{assumption 1} and \ref{assumption 2} we find $m_{k}$ and $\psi_{k}$ through the following sequence of choices:
  	\begin{enumerate}[label=(\roman*)]
  		\item $\theta_{k}\in C^{0}(X),\,\theta_{k}>0$ such that 
  		\begin{enumerate}[label=(\alph*), ref=(\roman{enumi}\alph*)]
  			\item  \label{1a}$\nn{\theta_{k}-\mathcal{D}(g_{k})}_{1}<\frac{\eta}{2^{k+2}\nn{\Delta_{k-1}}_{\infty}}$ 
  			\item  \label{1b}and $\nn{\theta_{k}-\mathcal{D}(g_{k})}_{1,P(m_{k-1},l)}<\min\,\{\frac{\epsilon'}{2},(\frac{\epsilon}{32})^{2}\}\nu_{k-1,l}$ for $l\in\lk_{m_{k-1}}$, 
  		\end{enumerate}
  		\item \label{3} $\lambda_{k}\in(0,\delta_{k})$ such that $\omega_{\theta_{k}}(2\lambda_{k})<\frac{\eta}{2^{k+2}\nn{\Delta_{k-1}}_{\infty}}$,
  		\item \label{mk3}and finally $m_{k}>m_{k-1}$ such that $\sup_{l\in\lk_{m_{k}}}\diam(P(m_{k},l))<\lambda_{k}$.
  	\end{enumerate}
    
    Define $\phi_{k}$ as the step function  
    $$
    \phi_{k}:=\sum_{l\in\lk_{m_{k}}}\left (\inf_{x\in P(m_{k},l)}\theta_{k}(x)\right)\ind_{P(m_{k},l)}.
    $$

  	We start the proof of Proposition with the following sublemma. 
    
    \begin{lemma}
    	\label{l: bounding by martingale}There is some function $\psi_{k}$ constant on $P(m_{k},l)$ for every $l\in\lk_{m_{k}}$ and satisfying 
    	\begin{enumerate}[label=(\alph*)]
    		\item \label{prop1} $\ln(\phi_{k})\leq -\epsilon'+\psi_{k}$ almost everywhere, 
    		\item \label{prop2}$\int_{P(m_{k-1},l)}\psi_{k}\dm=0$ for each $l\in\lk_{m_{k-1}}$,
    		\item \label{prop3} and $\nn{\psi_{k}}_{2}^{2}\leq 8$. 
    	\end{enumerate}
    \end{lemma}
    \begin{subproof}
    	 First of all, notice that for all $l\in\lk_{m_{k-1}}$ by (\ref{e: delta choice}), the assumption $g_{k}\in\V_{\delta_{k}}$, and \ref{assumption 1} we have 
    	 \begin{equation}
    	 	\label{e: derivative inequality}\ms{g_{k}(P(m_{k-1},l))}\leq (1+\min\{\frac{\epsilon'}{2},\frac{\epsilon}{4}\})\nu_{k-1,l}.
    	 \end{equation}
    	 On the other hand, since at least $\mu$-half of $P(m_{k-1},l)$ is covered by 
    	 sets of the form $P(m_{k},l')$ with $l'\in\lk'_{m}$ (item \ref{en4} of \ref{d: partitions}), it follows from condition \ref{assumption 2} in the choice of $g_{k}$ that: 
    	  \begin{equation}
    	  	 \label{inequality sufficiently small derivative} \nn{g_{k}}_{ac,P(m_{k-1},l)}\geq\frac{\epsilon\nu_{k-1,l}}{2}.
    	  \end{equation}
    	 Together with \eqref{e: derivative inequality}, this allows us to apply \cref{c: often enough negative enough} with $\frac{\epsilon}{2}$ in place of the $\epsilon$ in the statement of the Corollary and $\frac{\epsilon}{4}$ in place of $\delta$, and obtain the inequality
    	 $$\mu(\{x\in P(m_{k-1},l)\,|\,\mathcal{D}(g_{k})(x)\leq 1-\frac{\epsilon}{16}\})\geq\frac{\epsilon\nu_{k-1,l}}{16}.$$ 
  	   From the bound $\nn{\mathcal{D}(g_{k})-\theta_{k}}_{1,P(m_{k-},l)}\leq\frac{\epsilon'}{2}\nu_{k-1,l}$ in \ref{1b} it follows that
  	   \begin{equation}\label{inequality aproxximation}
  	   	\mu(\{x\in P(m_{k-1},l)\,|\,|\mathcal{D}(g_{k})(x)-\theta_{k}(x)|\geq \frac{\epsilon}{32}\})\leq \frac{\epsilon\nu_{k-1,l}}{32}.
  	   \end{equation}
  	   Write $A_{l}:=\{x\in P(m_{k-1},l)\,|\,\phi_{k}\leq 1-\frac{\epsilon}{32}\}$.  Note that this is a union of sets of the form $P(m_{k},l')$, $u\in\lk_{m_{k}}$, $P(m_{k},l')\subseteq P(m_{k-1},l)$. We can use inequalities \eqref{inequality sufficiently small derivative}  and \eqref{inequality aproxximation} to conclude
  	   \begin{align}
  	    	\label{e: small in a large enough set}
  	    	\begin{split}
  	    	 \mu(A_{l})\geq\mu(\{x\in P(m_{k-1},l)\,|\,\theta_{k}\leq 1-\frac{\epsilon}{32}\})\geq\frac{\epsilon\nu_{k-1,l}}{32}, 
  	    	\end{split}
  	    \end{align}
  	    since clearly $\mathcal{D}(x)\leq 1-\frac{\epsilon}{16}$ and $|\mathcal{D}(g_{k})(x)-\theta(x)|<\frac{\epsilon}{32}$ imply $\theta(x)<\frac{\epsilon}{32}$ for $x\in P(m_{k-1},l)$.
  	   
  	    Next, we show the following result $K_{l}:=\int_{P(m_{k-1},l)}\Lambda(\phi_{k})\dm$.
  	    \begin{lemma}
  	    	\label{bounding Kl}The inequalities $-\nu_{k-1,l}\leq K_{l}\leq\epsilon'\nu_{k-1,l}$ hold.
  	    \end{lemma}
        \begin{subsubproof}
        	We start by bounding from above a natural approximation of $K_{l}$ using \eqref{e: derivative inequality}:
          \begin{equation}
    	      \label{bound on int of lambda}\int_{P(m_{k-1},l)}\Lambda(\mathcal{D}(g_{k}))\dm \leq\int_{P(m_{k-1},l)}(\mathcal{D}(g_{k})-1)\dm=\frac{\epsilon'\nu_{k-1,l}}{2}. 
          \end{equation}
          Using the monotonicity of $\Lambda$ for the first inequality below, the fact that $\Lambda$ is $1$-Lipschitz for the third one, and \eqref{bound on int of lambda} together with item \ref{1b} for the last one we obtain: 
          \begin{equation}
              \label{e: kl}
              \begin{split}
                 K_{l}\leq&\int_{P(m_{k-1},l)}\Lambda(\theta_{k})\dm\leq\int_{P(m_{k-1},l)}\Lambda(\mathcal{D}(g_{k}))\dm+\nn{\Lambda(\mathcal{D}(g_{k}))-\Lambda(\theta_{k})}_{1,P(m_{k-1},l)}\leq\\ &\int_{P(m_{k-1},l)}\Lambda(\mathcal{D}(g_{k}))\dm+\nn{\mathcal{D}(g_{k})-\theta_{k}}_{1,P(m_{k-1},l)}\leq\frac{\epsilon'\nu_{k-1,l}}{2}+\frac{\epsilon'\nu_{k-1,l}}{2}=\epsilon'\nu_{k-1,l}.
              \end{split}
           \end{equation}
         The bound $K_{l}\geq-\nu_{k-1,l}$ follows trivially from the fact that $\Lambda\circ\phi>-1$ everywhere.  
        \end{subsubproof}

    Next, write
    \begin{equation*}
    	\gamma_{l}:=\max\left\{0, \frac{K_{l}+\epsilon'\nu_{k-1,l}}{\mu(A_{l})}\right \},\quad\quad \sigma_{l}:=\frac{\gamma_{l}\mu(A_{l})-K_{l}}{\nu_{k-1,l}}.
    \end{equation*}
    
    \begin{lemma}
    	\label{bound on gamma and sigma} For all $l\in\lk_{m_{k-1}}$ we have
    	\begin{itemize}
    		\item $\sigma_{l}\geq\epsilon'$,
    		\item $\{\sigma_{l},\sigma_{l}-\gamma_{l}\}\subseteq [-1,1]$. 
    	\end{itemize}
    \end{lemma}
    \begin{subsubproof}
    	Observe that for each $l\in\lk_{m_{k-1}}$ we have one of the two situations below:
    \begin{itemize}
    	\item $\gamma_{l}=\frac{K_{l}+\epsilon'\nu_{k-1,l}}{\mu(A_{l})}>0$ and $\sigma_{l}=\epsilon'$,
    	\item $\gamma_{l}=0$, the quantity $\frac{K_{l}+\epsilon'\nu_{k-1,l}}{\mu(A_{l})}$ is negative, and $\sigma_{l}=-\frac{K_{l}}{\nu_{k-1,l}}$.   
    \end{itemize}
  	
  	Using the upper bound in \cref{bounding Kl} in the first inequality below and inequality \eqref{e: small in a large enough set} in the second we get: 
    \begin{equation*}
    	\gamma_{l}\leq\frac{2\epsilon'\nu_{k-1,l}}{\mu(A_{l})}\leq\frac{ 64\epsilon'}{\epsilon}=\hat{\epsilon}.
    \end{equation*}
    	If $\gamma_{l}>0$, then $\sigma_{l}=\epsilon'$ and the inclusion $\{\sigma_{l},\sigma_{l}-\gamma_{l}\}\subseteq[-1,1]$ is clear (for $\epsilon$ small enough).
    	In the second case we have 
    	 \begin{equation*}
    	 	\epsilon '\leq \sigma_{l}=\sigma_{l}-\gamma_{l}=-\frac{K_{l}}{\nu_{k-1,l}}\leq 1,
    	 \end{equation*}
    	where the first inequality follows from $-K_{l}+\epsilon'\nu_{k-1,l}<0$ and the second from \cref{bounding Kl}.
    	
    \end{subsubproof}
   To conclude the proof of \cref{l: bounding by martingale}, define
   $$
     \psi_{k}:=\Lambda(\phi_{k})+\sum_{l\in\lk_{m_{k-1}}}(\sigma_{l}\ind_{P(m_{k-1},l)}-\gamma_{l}\ind_{A_{l}}).
   $$
   We proceed to check that $\psi_{k}$ satisfies the properties in the statement of \cref{l: bounding by martingale}, thereby concluding its proof. 
   \begin{itemize}
   	\item[\textbf{Item \ref{prop1}}.] The inequality 
   	\begin{equation*}
   		\ln(\phi_{k})\leq-\epsilon'+\psi_{k}
   	\end{equation*}
   	 holds almost everywhere. Indeed, by definition, 
   	$\hat{\epsilon}=-\frac{\epsilon}{32}-\ln(1-\frac{32}{\epsilon})$, that is $\hat{\epsilon}=\Gamma(1-\frac{\epsilon}{32})$, where
   	\begin{equation*}
   		 \Gamma(y)=(y-1)-\ln(y).
   	\end{equation*}
   	It is not hard to see that $\Gamma$ is strictly decreasing on $(0,1]$. 
   	For all for all $x\in A_{l}$ we have $\phi(x)\leq\theta(x)\leq 1-\frac{\epsilon}{32}$ and $\Lambda(\phi(x))=x-1$. Therefore for all $x\in A_{l}$ outside a set of measure zero we have: 
   	\begin{equation*}
   		\ln(\phi_{k}(x))\leq\Lambda(\phi(x))-\hat{\epsilon}\leq \Lambda(\phi(x))-\gamma_{l}=-\sigma_{l}+\psi(x)\leq-\epsilon'+\psi(x).
   	\end{equation*}
   	The inequality between first and last terms can be proven for almost all $x\in P(m_{k-1},l)\setminus A_{l}$ in a similar way. Recall that $\epsilon'\leq\sigma_{l}$ is from \cref{bound on gamma and sigma}.
   	\item[\textbf{Item \ref{prop2}}.] For all $l\in\lk_{m_{k-1}}$ the choice of $\sigma_{l}$ 
   	implies that 
   	\begin{equation*}
   		\int_{P(m_{k-1},l)}\psi_{k}\dm=K_{l}+\sigma_{l}\nu_{k-1,l}-\gamma_{l}\mu(A_{l})=0.
   	\end{equation*}
   	\item[\textbf{Item \ref{prop3}}.]  Finally, the bound $\nn{\psi_{k}}_{2}\leq 8$ holds, as a consequence of \cref{o: small variance after perturbation}, which applies with 
   	$
   		\rho:=\sum_{l\in\lk_{m_{k-1}},l}(\sigma_{l}\ind_{P(m_{k-1},l)}-\gamma_{l}\ind_{A_{l}}) 
   	$
   	 using the conclusion of \cref{bound on gamma and sigma} and inequality $\nn{\phi_{k}}_{1,P}\leq 2$.
   \end{itemize} 
    \end{subproof}
   	
   	\subsection*{End of the proof of \cref{p: almost independent}: $L^{1}$-approximate bound of the derivative.} We are now in a position to finish the proof of \cref{p: almost independent}. Assume that the process has run for $n$ steps, producing elements $g_{1},\dots g_{n}\in \gg$, together with all the other data. Write $\hat{g}_{k}:=g_{k+1}\circ\dots g_{n}$ for all $0\leq k\leq n$. First of all, set 
  	\begin{alignat*}{2}
  		&\Theta_{k}:=\Delta_{k-1}(\mathcal{D}(g_{k}\circ\hat{g}_{k})-\phi_{k}\mathcal{D}(\hat{g}_{k}))\in L^{1}(X,\mu),\quad && 1\leq k\leq n,\\
  		&\Xi_{k}:=\Delta_{k}\mathcal{D}(\hat{g}_{k}), \quad && 0\leq k\leq n.
  	\end{alignat*} 
  	For all $1\leq k\leq n$ noting that $g_{k}\circ\hat{g}_{k}=\hat{g}_{k-1}$ we have the validity almost everywhere of
    \begin{equation}
    	\label{step inequality}\Xi_{k-1}-\Theta_{k}=\Delta_{k-1}\phi_{k}\mathcal{D}(\hat{g}_{k})\leq \Delta_{k-1} e^{-\epsilon'+\psi_{k}}\mathcal{D}(\hat{g}_{k})=\Delta_{k}\mathcal{D}(\hat{g}_{k})=\Xi_{k},
    \end{equation}
    where the first equality follows from some simple algebraic manipulation, the inequality follows from the conclusion of Lemma \ref{l: bounding by martingale} and the final equality from the definition $\Delta_{j}:=e^{-j\epsilon'+\sum_{l=1}^{j}\psi_{l}}$.
    
  	Adding inequality \eqref{step inequality} for all values $1\leq k\leq n$ we get the validity almost everywhere of  
  	$$
  	\mathcal{D}(g_{1}\circ\dots g_{n})=\Xi_{0}\leq \Theta+\Xi_{n}=\Theta+e^{-n\epsilon'+\sum_{k=1}^{n}\psi_{l}},
  	$$
  	where $\Theta:=\sum_{k=1}^{n}\Theta_{k}$. To conclude the proof of \cref{p: almost independent} it suffices to show that $\nn{\Theta}_{1}\leq\eta$, which boils down to the following:
  	\begin{lemma}
  		\label{inequality right side}The inequality $\nn{\mathcal{D}(g_{k}\circ\hat{g}_{k})-\phi_{k}\mathcal{D}(\hat{g}_{k})}_{1}\leq\frac{\eta}{2^{k}\nn{\Delta_{k-1}}_{\infty}}$ holds for $1\leq k\leq n$. 
  	\end{lemma}
     \begin{subproof}
     	 We just need to verify that we may apply \cref{o: approximation use} with $\frac{\eta}{2^{k+1}\nn{\Delta_{k-1}}_{\infty}}$  in place of the values of the $\epsilon$ and $\epsilon'$ in the statement of \ref{o: approximation use}. Recall that by construction $\sup_{l\in\lk_{m_{k}}}\diam(P(m_{k},l))<\lambda_{k}$, and $\phi_{k}$ agrees with $\inf_{x\in P(m_{k},l)}\theta_{k}(x)$ on each of those sets. 
     	  
  	   The $L^{1}$ bound in the assumptions of \cref{o: approximation use} follows from \ref{1a}, combined from the $L^{\infty}$ bound resulting from the choice of $\lambda_{k}$ in \ref{3}, as well as the properties of $\phi_{k}$ mentioned above:
  	   \begin{equation*}
  	   	\nn{\mathcal{D}(g_{k})-\phi_{k}}_{1}\leq \nn{\mathcal{D}(g_{k})-\theta_{k}}_{1}+\nn{\theta_{k}-\phi_{k}}_{\infty}\leq 2\frac{\eta}{2^{k+2}\nn{\Delta_{k-1}}_{\infty}}=\frac{\eta}{2^{k+1}\nn{\Delta_{k-1}}_{\infty}}.
  	   \end{equation*}
  	   
  	   In order to verify the conditions pertaining to the $\epsilon'$ in \cref{o: approximation use}, observe first that for all $1\leq s\leq n$ we have $\delta_{s+1}\leq\frac{\lambda_{s}}{2}<\frac{\delta_{s}}{2}$, as well as $g_{s}\in\V_{\delta_{s}}$. It follows that $\hat{g}_{k}\in\V_{2\delta_{k+1}}\subseteq\V_{\lambda_{k}}$. 
    \end{subproof}
  	From \cref{inequality right side} we get:
     	\begin{align}
  	\label{e: first l1 bound}
  	\begin{split}
  		\nn{\Theta_{k}}_{1}\leq&\nn{\Delta_{k-1}}_{\infty}\nn{((\mathcal{D}(g_{k})\circ\hat{g}_{k}-\phi_{k})\mathcal{D}(\hat{g}_{k})}_{1}
  		\leq \nn{\Delta_{k-1}}_{\infty}\frac{2\eta}{2^{k+1}\nn{\Delta_{k-1}}_{\infty}}=\frac{\eta}{2^{k}},
    \end{split}
  	\end{align}
  	from which it is immediate that $\nn{\Theta}_{1}<\eta$, concluding the proof of \cref{p: almost independent}.
  \end{proof}

	\begin{lemma}
	  	\label{c: almost final lemma} Assume that $(X,\mu)$ is well-covered and let $\tau$ be a group topology on $\hh(X,\mu)$ satisfying $\tau_{co}\subsetneq \tau$ and $\tau_{ac}\nsubseteq \tau$. Then for any $\V\in\idn$ and any $\delta>0$ there exists some $g\in\V$ such that 
	  	$\mu(\{x\in X\,|\,\mathcal{D}(g)(x)\geq\delta\})<\delta$. 
	\end{lemma} 
  \begin{proof}
  	Pick $\epsilon_{0}>0$ such that for every $\V\in\idn$ there is some $g\in\V$ with $\nn{g}_{ac}>\epsilon_{0}$ and write $\epsilon=\frac{\epsilon_{0}}{64}$ and $\epsilon':=\frac{\epsilon}{64}(\frac{\epsilon}{32}-\ln(1-\frac{\epsilon}{32})-1)$.
  	Let also $(\mathcal{P}_{m},\lk_{m},\lk'_{m})_{m\geq 0}$ be an $\frac{\epsilon_{0}}{4}$-filling good covering sequence for $(X,\mu)$. 
  	Let 
  	$$\nu:=\min\left\{\frac{2}{\epsilon'},\frac{\delta}{2}\right\},\quad\quad\quad N_{0}=\max\,\left\{\,\bigg\lceil\frac{2\ln(\sfrac{2}{\delta})}{\epsilon'}\bigg\rceil,\,N(8,\nu)\,\right\},$$
  	where $N :\R_{>0}^{2}\to\mathbb{N}$ is the function given by Corollary \ref{c: Chebyshev martingale}. Pick also $\mathcal{U}=\mathcal{U}^{-1}\in\idn$ with $\mathcal{U}^{N_{0}}\subseteq\V$.
  	Using Lemma \ref{c: uniform derivative} and the assumption that $\tau_{co}\subseteq\tau$ we can inductively choose $g_{1},\dots g_{N_{0}}\in\mathcal{U}$ 
  	such that for some $0<\eta<(\frac{\delta}{2})^{2}$ and 
  	$\delta_{k}:=\df_{\eta,k}(g_{1},\dots g_{k-1})$, $m_{k}=\nf_{\eta,k}(g_{1},\dots g_{k})$, 
   where $\df_{\eta,k},\nf_{\eta,k}$ are given by Proposition \ref{p: almost independent}, then the following conditions are satisfied: 
   	\begin{itemize}
   		\item $g_{k}\in\V_{\delta_{k}}$,
   		\item $\nn{g_{k}}_{ac,P(m_{k},l)}>\epsilon\mu(P(m_{k},l))$ for $l\in\lk'_{m_{k}}$. 
   	\end{itemize}
  
  	If we write $g:=g_{1}\circ\cdots g_{n}$, 
  	then \ref{p: almost independent} yields $\Theta\in L^{1}(X.\mu)$ with $\nn{\Theta}_{1}<(\frac{\delta}{2})^{2}$ and a sequence of martingale increments
  	$\psi_{1},\dots \psi_{n}$ such that $\nn{\psi_{l}}_{2}^{2}\leq 8$ and 
  	$\mathcal{D}(g)\leq \Theta+e^{-N_{0}\epsilon'+\sum_{l=1}^{N_{0}}\psi_{l}}$ almost everywhere. Using one of the bounds on $N_{0}$ in each of the last two inequalities below, we get: 
  	\begin{align*}
  		&\mu(\{x\in X\,|\,(\mathcal{D}(g)-\Theta)(x)\geq\frac{\delta}{2}\})\leq\mu(\{x\in X\,|\,e^{-N_{0}\epsilon'+\sum_{l=1}^{N_{0}}\psi_{l}}\geq\frac{\delta}{2}\})\leq\\
  		&\mu(\{x\in X\,|\,-\epsilon'+\frac{1}{N_{0}}\sum_{l=1}^{N_{0}}\psi_{l}(x)\geq\frac{\ln(\frac{\delta}{2})}{N_{0}}\})\leq
  		\mu(\{x\in X\,|\,|\frac{1}{N_{0}}\sum_{l=1}^{N_{0}}\psi_{l}(x)|\geq\frac{\epsilon'}{2}\})\leq\frac{\delta}{2}   
  	\end{align*}
  	On the other hand, $\nn{\Theta}_{1}<(\frac{\delta}{2})^{2}$ implies that 
  	$\mu(\{x\in X\,|\,\Theta(x)>\frac{\delta}{2}\})\leq\frac{\delta}{2}$ and thus
  	$$\mu(\{x\in X\,|\,\mathcal{D}(g)\geq\delta\})\leq\mu(\{x\in X\,|\,(\mathcal{D}(g)-\Theta)(x)\geq\frac{\delta}{2}\})+\mu(\{x\in X\,|\,\Theta(x)\geq\frac{\delta}{2}\})\leq\delta$$ 
  \end{proof}
  
  Using the regularity of the measure $\mu$ together with Observation \ref{o: change of variable} one easily concludes:
  \begin{corollary}
  	\label{c: shrinking sets}Under the assumptions of Lemma \ref{c: almost final lemma} for any $\V\in\idn$ and any $\delta>0$ there is some open set $V\subseteq X$ and some $g\in\V$ such that 
  	$\mu(V)\geq 1-\delta$ and $\mu(g(V))\leq\delta$. 
  \end{corollary}

  \section{Separable topologies on $\hh(X,\mu)$ in the Oxtoby-Ulam case} 
  
  \label{s: separability}
  Throughout this section we let $X$ be a compact manifold and $\mu$ an OU measure on $X$ and $\gg:=\hh(X,\mu)$. Without loss of generality we assume that $\mu$ is a probability measure.
  
     \begin{definition}
     	  \label{d: nice collections^}Assume that $X\neq[0,1]$. Given $\epsilon>0$ by an $\epsilon$-fine collection we intend a finite disjoint collection $\{D_{l}\}_{l\in\Lambda}$ of embedded balls and boundary half-balls such that: 
       \begin{itemize}
        	\item $\diam(D_{l})\leq\epsilon$ for for all $l\in\Lambda$,
        	\item and $\mu(\bigcup_{l\in\Lambda}D_{l})\leq\epsilon$. 
       \end{itemize}
       By an $\epsilon$-fine open subset of $X$ we intend one of the form $\bigcup_{l\in\Lambda}\mathring{D}_{l}$ for some $\epsilon$-fine collection $\{D_{l}\}_{l\in\Lambda}$.
     \end{definition}
  
   The goal of this section is to prove the result below, as well as a variation of it for $X=[0,1]$. 
  \begin{lemma}
   	\label{l: separable topologies ou case} 
   	If $X\neq[0,1]$, then there is some number $N$ such that for any separable group topology $\tau$ on $\gg:=\hh(X,\mu)$ and every $\V\in\idn$ there exists some $\epsilon>0$ such that 
   	for any $\epsilon$-fine open subset $V$ of $X$ we have $\gg_{V}\subseteq\mathcal{V}$.  
  \end{lemma}
  \begin{proof}

   By continuity of multiplication, it suffices to show, after fixing $\epsilon$, that for any arbitrary $\mathcal{U}=\mathcal{U}^{-1}\in\idn$ the result is true with $\V$ replaced by $\mathcal{U}^{14}$ in the conclusion. Fix such $\mathcal{U}$.
   For the following see Lemma 3.2 in \cite{mann2016automatic}, which uses only the Baire property.
   \begin{lemma}
   	\label{l: approximating neighborhoods}The closure of $\mathcal{U}$ contains the standard neighborhood of the identity $\W_{\delta}$ of $\tau_{ac}$ (Definition \ref{d: neighborhoods}) for some $\delta>0$. 
   \end{lemma}
   Take $\epsilon:=\frac{\delta}{4}$ for the above $\delta$ and fix an $\epsilon$-fine collection $\{D_{l}\}_{l\in\Lambda}$. The following follows from Observations \ref{o: pasting} by the same argument as in Lemma 3.8 in \cite{mann2016automatic}.
   \begin{lemma}
   	\label{l: extensions in neighborhood} Assume that for all $l\in\Lambda$ we are given a countable disjoint collection $\{D_{l,k}\}_{k\geq 0}$ of sets of the same type as $D_{l}$. Then there is some choice function
   	$\rho:\Lambda\to\N$ such that for all $l\in\Lambda$ and any
   	choice of elements $g_{l}\in\gg_{D_{l,\rho(l)}}$ there is some $g\in\gg_{\bigcup_{l\in\Lambda,k\geq 0}D_{l,k}}\cap\mathcal{U}^{2}$ such that $g_{\restriction D_{l,\rho(l)}}$ for all $l$.   
   \end{lemma}
    
   Using Corollary \ref{c: uniform perfectness} and Lemma \ref{l: translation} together with the argument in the proof of Lemma 
   3.7 in \cite{mann2016automatic} and the inclusion $\tau_{co}\subseteq\tau_{ac}$ we obtain the following refinement of Lemma \ref{l: approximating neighborhoods}: 
   \begin{lemma}
   	 \label{l: extensions 2} Let $D_{l,k}\subseteq X$, $\lambda\in\Lambda, k\geq 0$  be as in Lemma \ref{l: extensions in neighborhood} and
    	$\rho:\Lambda\to\N$ the choice function provided therein. Assume that for all $l\in\Lambda$ a set $D_{l}'\subseteq \mathring{D}_{l,\rho(l)}$ of the same type as $D_{l}$ is given. Then for any	choice of elements $g_{l}\in\gg_{\mathring{D}_{l,\rho(l)}}$ there is $g\in\gg_{\bigcup_{l\in\Lambda}D'_{l}}\cap\mathcal{U}^{12}$ 
   	such that $g_{\restriction D'_{l,\rho(l)}}=g_{l}$ for all $1\leq l\leq r$.
   \end{lemma}

   To conclude the proof of \ref{l: separable topologies ou case} let $\{D'_{l}\}_{l\in\Lambda}$ be as given by 
   Lemma \ref{l: extensions 2}. We can use Observation \ref{o: approximation use} to find $h_{0}\in\W_{\delta}$ such that 
   $h_{0}(D_{l})\subseteq\mathring{D'_{l}}$ for all $l\in\Lambda$ (recall that $\Lambda$ is finite) and by Lemma \ref{l: approximating neighborhoods} (in fact, we only need $\tau_{co}$-density) we may assume that $h_{0}\in\mathcal{U}$. It follows that $\gg_{\bigcup_{l\in\Lambda}D_{l}}\subseteq\mathcal{U}^{h_{0}}\subseteq\mathcal{U}^{14}$.
  \end{proof}
   
     \begin{definition}
     	  \label{d: nice collections 2}Let $\mu$ be the Lebesgue measure on $X=[0,1]$, $\epsilon>0$ and $\omega:[0,1]\to\R_{\geq 0}$ some continuous function with 
     	  $\omega^{-1}(0)=\{0,1\}$. By an $(\omega,\epsilon)$-nice collection we intend a countable disjoint collection $\{D_{l}\}_{l\in\Lambda}$ of closed subintervals of $(0,1)$ with: 
       \begin{itemize}
        	\item $\diam(D_{l})\leq sup_{x\in D_{l}}\omega(x)$ for all $l\in\Lambda$,
        	\item and $\mu(\bigcup_{l\in\Lambda}D_{l})\leq\epsilon$. 
       \end{itemize}
       We will denote any set of the form $\bigcup_{l\in \Lambda}\mathring{D_{l}}$ for some collection $\{D_{l}\}_{l\in\Lambda_{l}}$ as above as an $(\omega,\epsilon)$-nice set. 
      \end{definition}
   
  \begin{lemma}
   	\label{l: separable topologies ou case 2} 
   	Let $X=[0,1]$ and $\mu$ be the Lebesgue measure. Then there is some number $N$ such that for any separable group topology $\tau$ on $\gg:=\hh(X,\mu)$ and every $\V\in\idn$ there are $\omega,\epsilon$ as in Definition \ref{d: nice collections 2} such that 
   	for any $\epsilon$-fine $V\subseteq X$ we have $\gg_{V}\subseteq\mathcal{V}$.  
  \end{lemma}
  \begin{proof}
  	Let $\tilde{d}$ be some compatible proper metric on $(0,1)$. 
  	One can define a modified metric $\tilde{\rho}$ on $\hh(0,1)$ by
  	$$
  	\tilde{\rho}(g,h):=\min\{\,1,\,\sup_{x\in(0,1)}\tilde{d}(g(x),h(x))+\sup_{x\in(0,1)}\tilde{d}(g^{-1}(x),h^{-1}(x))+\nn{\mathcal{D}(g)-\mathcal{D}(h)}_{1}+\nn{\mathcal{D}(g^{-1})-\mathcal{D}(h^{-1})}_{1}\}.
  	$$ 
  	The same argument as in Proposition \ref{p: generalization} shows that $\tilde{\rho}$ is complete and induces a group topology. It 
  	is not separable, but the resulting topological space will have the Baire property nevertheless, so that \ref{l: approximating neighborhoods} will hold with the topology induced by $\tilde{\rho}$ in place of $\tau_{ac}$. From this point on, the argument continues in the same fashion as in Lemma \ref{l: separable topologies ou case}. Alternatively, one may use the existence of generic conjugacy classes \cite{ihli2022generation} and the argument in the proof of automatic continuity for $\Homeo^{+}(\R)$ found in \cite{rosendal2006automatic}. 
  \end{proof}

  \section{Proof of Theorem \ref{t: main}.}
  \label{s:proof of b}

   Given an open set $V\subseteq X$ we let $\cdiam(V)$ stand for the supremum of the diameter of all connected components of $V$. 
  \begin{proposition}
  	\label{p: first technical statement} Let $(X,\mu)$ be as in Proposition \ref{p: generalization} and let $\gg:=\hh(X,\mu)$ and additionally assume that
  	\begin{enumerate}[label=(\alph*)]
  		\item \label{conda}$(X,\mu)$ is well-covered,
  		\item \label{condb} $\gg$ has the local fragmentation property for $\tau_{co}$,
  		\item \label{condc} we are given some collection $\mathscr{C}$ of open sets in $X$ and $N>0$ with the following properties:
  		 \begin{enumerate}[label=(\roman*), ref=(\alph{enumi}\roman*)]
  		 	\item \label{condition i}for any $U\in\mathscr{C}$ and $\mathcal{F}\subseteq\pi_{0}(U)$ we have $\bigcup\mathcal{F}\in\mathscr{C}$, 
  		 	\item \label{condition ii}for any $\epsilon>0$ there is a covering of $X$ by sets  $V_{1},\dots V_{N}\in\mathscr{C}$ with $\cdiam(V_{j})\leq\epsilon$.
  		 \end{enumerate}
  	\end{enumerate}
     Then $\tau_{co}$ is the only group topology $\tau$ on $\gg$ satisfying
     \begin{enumerate}[label=(\arabic*)]
     	\item \label{propa}$\tau_{ac}\nsubseteq\tau$, $\tau_{co}\subseteq\tau$ 
     	\item \label{propc} For any $\V\in\idn$ there is $\epsilon>0$ such that for any open $V\in\mathscr{C}$ with $\bar{\mu(V)}\leq\epsilon$ and $\cdiam(V)\leq\epsilon$ we have $\gg_{V_{j}}\subseteq\V$. 
     \end{enumerate}  
     In particular, there is no group topology $\tau$ on $\hh(X,\mu)$ with $\tau_{co}\subsetneq\tau\subsetneq\tau_{ac}$.
  \end{proposition}
  \begin{proof}
  	Assume $\tau_{ac}\nsubseteq\tau$ and $\tau_{co}\subseteq\tau$ and pick some arbitrary $\V\in\idn$ 
  	and $\mathcal{U}=\mathcal{U}^{-1}\in\idn$ with $\mathcal{U}^{2N+2}\subseteq\V$. Let $\epsilon>0$ be the constant associated with $\mathcal{U}$ by \ref{propc}. By Corollary \ref{c: shrinking sets}
  	there is some open set $V\subseteq X$ and some $g\in\mathcal{U}$ such that $\mu(V)\geq 1-\frac{\epsilon}{2}$ and $\mu(g(V))\leq\frac{\epsilon}{2}$. Using the regularity and continuity of the measure one can easily find open sets $V',V''$ with 
  	$\bar{V'}\subseteq V$, $V'\cup V''=X$ and $\mu(\bar{V'})\leq\epsilon$.  
  	
  	Now, fix some $\delta>0$ such that $\delta\leq\epsilon$, $\omega_{g}(\delta)\leq\epsilon$ and $\delta$ is smaller than the Lebesgue number of the open cover $\{V',V''\}$. Choose $V_{1},\dots V_{n}\in\mathscr{C}$ with $X=\bigcup_{j=1}^{N}$ and $\cdiam(V_{j})<\delta$. For each $1\leq j\leq N$ define $V'_{j}\in\mathscr{C}$ as the union of the connected components of $V_{j}$ contained in $V'$ and define $V''_{j}\in\mathscr{C}$ analogously. Our choice of $\delta$ ensures that $V_{j}=V'_{j}\cup V''_{j}$
  	and that $\cdiam(g(V'_{j}))\leq\epsilon$. The choice of $\epsilon$ implies that $\gg^{c}_{g(V'_{j})}\cup\gg^{c}_{V''_{j}}\subseteq\mathcal{U}$.  
  	
  	The fragmentation property implies the existence of some $\eta>0$ such that
  	$$
  	\mathcal{V}_{\eta}\subseteq\gg^{c}_{V''_{1}}\cdots\gg^{c}_{V''_{N}}\gg^{c}_{V'_{1}}\cdots\gg^{c}_{V'_{N}}\subseteq\mathcal{U}^{N}(\mathcal{U}^{g})^{N}\subset\mathcal{U}^{2N+2}\subseteq\V,
  	$$
    and so we are done. 
  \end{proof}

  Our main theorem follows. 
  \begin{proof}[Proof of Theorem \ref{t: main}]
   Assume first $(X,\mu)$ is an OU pair with $X\neq[0,1]$, $\gg:=\hh(X,\mu)$. 
   We need to check that conditions \ref{conda}, \ref{condb} of the proposition are satisfied by $(X,\mu)$ and that the separability hypothesis on $\tau$ implies condition \ref{propc} on $\tau$, where $\mathscr{C}$ is the collection of all finite unions of balls and boundary half-balls (note that \ref{propc} is trivial in case $\tau\subseteq\tau_{ac}$ by Lemma \ref{l: l1bound} for the collection of all open sets). This follows from Corollary \ref{c: ou filling} (for filling partitions), Lemma \ref{l: fragmentation} (for local fragmentation) and \cref{l: nice colorable coverings} in combination with Lemma \ref{l: separable topologies ou case} (for property \ref{propc}).
   
   The case \ref{Cantor case} ($X=\mathcal{C}$) of the statement the first two properties are trivial from the assumption, the equivalent of Fact \ref{f: colorable coverings} is clear as well, and the equivalent of Lemma \ref{l: separable topologies ou case} can be shown using an simplified version of the same argument much as in \cite{rosendal2006automatic} (using the commutator trick instead of the existence of generic conjugacy classes).  
   
   The case $X=[0,1]$ is somewhat more complicated. One can proceed along the lines of Theorem \ref{p: first technical statement} taking as $\mathscr{C}$ the collection of families of disjoint intervals and requiring in \ref{propc} that the connected components of $V$ are $(\omega,\epsilon)$ for some $\omega,\epsilon$ as in Definition \ref{d: nice collections 2}. That this variant of \ref{propc} is satisfied is the content of Lemma \ref{l: separable topologies ou case 2}. Here we need to use Corollary \ref{c: fragmentation finer} instead of Lemma \ref{l: fragmentation}, which yields a weaker conclusion, but then the argument involving infinite combinatorics found in \cite{rosendal2006automatic}, p.13 can be applied to reach the desired conclusion.
    
  \end{proof}

  \section{Roelcke precompactness: proof of Theorem \ref{t: main2} and Theorem \ref{t: main3}}
  \label{s:roelcke}

 Let us now show Theorem \ref{t: main2}.
   \begin{customthm}{\ref{t: main2}}
  	\mainthmbodyy
  \end{customthm}
   \begin{proof}
 For simplicity, assume that $\mu$ is the Lebesgue measure and $d$ the standard distance. 	
   \begin{observation}
  	\label{o: key observation} Suppose that we are given $u,\tilde{u},g,g'\in\gg$ with $g=\tilde{u}g'u$ and $g'\in\V_{\delta}$, as well as some finite collection of disjoint open intervals $\{J_{i}\}_{1\leq i\leq k}$ in $X$. Then
  	$$
  	 \nn{u}_{ac}+\nn{\tilde{u}}_{ac}\geq\sum_{i=1}^{k}|\ms{J_{i}}-\ms{g(J_{i})}|-2k\delta. 
  	$$
  \end{observation}
  \begin{subproof}
    For $1\leq i\leq k$ let $\Delta_{i}:=|\ms{J_{i}}-\ms{u(J_{i})}|$. Let also $J'_{i}:=g'u(J_{i})$ and $\Delta'_{i}:=|\ms{J'_{i}}-\ms{\tilde{u}(J'_{i})}|$.	
    Since $g'\in\V_{\delta}$, it is clear that $|\ms{J_{i}}-\ms{g(J_{i})}|\leq 2\delta+\Delta_{i}+\Delta'_{i}$ for all $i$. But then: 
    $$
    \sum_{i=1}^{k}|\ms{J_{i}}-\ms{g(J_{i})}|\leq 2k\delta+\sum_{i=1}^{k}\Delta_{i}+\sum_{i=1}^{k}\Delta'_{i}\leq 2k\delta+\nn{u}_{ac}+\nn{\tilde{u}}_{ac},
    $$
    as needed. 
  \end{subproof}

   	Let $(g_{n})_{n\in\N}\subseteq\gg$ a sequence such that $g_{n}\in\V_{\delta_{n}}$ for some sequence $(\delta_{n})_{n}\subseteq\R_{>0}$ converging to zero, but such that $(g_{n})_{n}$ does not converge to $\Id$ in $\tau_{ac}$. Up to passing to a subsequence, we may assume $\nn{g_{n}}\geq\epsilon$ for some $\epsilon>0$. If the conclusion of the theorem does not hold, then up to replacing $(g_{n})_{n}$ by a subsequence we may assume that for any $n>0$ there are $u_{n},\tilde{u}_{n}\in\V_{\frac{\epsilon}{8}}$ such that  
   	$g_{n}=u_{n}g_{0}\tilde{u}_{n}$. An approximation argument such as in Observation \ref{o: approximation use} yields the existence of some finite disjoint collection of intervals $J_{1},\dots J_{k}\subseteq X$ with disjoint interiors such that $\sum_{i=1}^{k}|\ms{g_{0}(J_{i})}-\ms{J_{i}}|\geq\frac{\epsilon}{2}$, which together with Observation \ref{o: key observation} leads to a contradiction:   	
   	$$
   	\frac{\epsilon}{4}\geq\nn{u_{n}}_{ac}+\nn{\tilde{u}_{n}}_{ac}\geq \frac{\epsilon}{2}-2k\delta_{n}\underset{n\to\infty}{\loongrightarrow} \frac{\epsilon}{2}.
   	$$       
    \end{proof}
   
   Key to making the proof above work is the fact that the condition $g\in\V_{\delta}$ implies some uniformly bounded distortion under the action of $g$ for all intervals. The absence of a similar property for the Cantor space with any reasonably homogeneous Borel probability measure is what makes Theorem \ref{t: main3} possible. Before proving it, we state two simple preliminary results.  
   
   \begin{fact}
   	\label{f: approximation} Let $\mu$ be a Borel probability measure on the Cantor space $\mathcal{C}$, $M>0$ a positive integer and  $\phi\in L^{1}(\mathcal{C},\mu)$, $\phi>0$ such that $\nn{\phi-\phi\cf_{\phi^{-1}(0,1]}}_{1}\leq\frac{1}{2M}$. Then 
   	there is a step function 
   	$\psi:=\sum_{j=1}^{3M}\frac{k}{3M}\cf_{A_{j}}$ where the $A_{j}$ are clopen and $\nn{\phi-\psi}_{1}\leq\frac{1}{M}$. 
   \end{fact}

   \begin{observation}
   	\label{o: homoteties} Let $K$ be a countable subfield of $\mathbb{R}$ and $\gg=\hh(\mathcal{C},\mu_{K})$. Then for any two partitions $\{A_{i}\}_{i=1}^{r}$ and $\{B_{i}\}_{i=1}^{r}$ of $\mathcal{C}$ into clopen sets there exists some 
   	$g\in \gg$ such that $g(A_{i})=B_{i}$ for all $1\leq i\leq r$ and $\mathcal{D}(g)_{\restriction A_{i}}$ is constantly equal to 
   	$\frac{\mu_{k}(B_{i})}{\mu_{K}(A_{i})}$. 
   \end{observation}
   \begin{proof}
   	 	The assumption that $K$ is a field, together with the uniqueness of Fr\"aiss\'e limits, implies that for any clopen set $A\subseteq\mathcal{C}$ there is a homeomorphism  
   	 	$h:\mathcal{C}\simeq A$ such that $h_{*}\mu_{K}=\mu_{K}(A)(\mu_{K})_{\restriction A}$, from which the result easily follows (compare with Observation \ref{o: ball correction for any}). 
   \end{proof}

   \newcommand{\nil}[0]{\emptyset}
   \begin{customthm}{\ref{t: main3}}
   	\mainthmbodyyy
   \end{customthm}
   \begin{proof}
   	\newcommand{\set}[0]{\{1,\dots r\}}
   	\newcommand{\cl}[0]{\mathbb{B}}
   	    \newcommand{\A}[0]{\hat{A}}
    \newcommand{\B}[0]{\hat{B}}
    \renewcommand{\aa}[0]{\A^{0,\pm}_{i,j}(k)}
    \newcommand{\bb}[0]{\B^{0,\pm}_{i,j}(k)} 
    \renewcommand{\AA}[0]{\A^{1,\pm}_{i,j}(k)}
    \newcommand{\BB}[0]{\B^{1,\pm}_{i,j}(k)} 
 
   	Let $\cl$ denote the boolean algebra of clopen subsets of $\mathcal{C}$, $\gg:=\hh(\mathcal{C},\mu)$. To prove $(\gg,\tau_{ac})$ is Roelcke precompact, it suffices to show that for all $\W\in\mathcal{N}_{\tau_{ac}}(\Id)$ 
   	and any infinite sequence $(g_{n})_{n\geq 0}\subseteq\gg$ there is some infinite subsequence $(g_{n_{k}})_{k\geq 0}$ such that $g_{n_{l}}\in \W g_{n_{k}}\W$ for all $k\neq l$. Without loss of generality there is some integer constant $N>0$ and some partition $\{A_{j}\}_{j=1}^{r}$ of $\mathcal{C}$ into clopen sets such that writing $\mh$ for the subgroup consisting of all the elements of $\gg$ preserving each $A_{j}$ setwise, we have that $\W$ consists precisely of those elements $g\in\mh$ satisfying
   	$\nn{g}_{ac}\leq\frac{1}{N}$. For $m\geq 0$, $i,j\in\{1,\dots r\}$ let
   	\begin{gather*}
   		A^{m}_{i,j}:=A_{i}\cap g_{m}^{-1}(A_{j})\quad\quad
   			A^{m,+}_{i,j}:=\{x\in A_{i,j}^{m}\,|\,\mathcal{D}(g_{m})(x)\geq 1\} \quad\quad A^{m,-}_{i,j}:=\{x\in A_{i,j}^{m}\,|\,\mathcal{D}(g_{m})(x)< 1\}
   	\end{gather*}

   	\newcommand{\C}[0]{2K}
   	\newcommand{\K}[0]{3K}
    Write $K:=6Nr^{2}$. For $m\geq 0$ and $(i,j)\in\{1,\dots r\}^{2}$ we can choose a partition $\{\hat{A}^{m,+}_{i,j},\hat{A}^{m,-}_{i,j}\}\subseteq\mathbb{B}$ of $A_{i,j}^{m}$ such that if we write $D_{i,j}^{m,+}:=A_{i,j}^{m,+}\triangle\hat{A}_{i,j}^{+}=A_{i,j}^{m,-}\triangle\hat{A}_{i,j}^{-}$, then 
   	
   	\begin{equation}
   		\label{e: good approx}\nn{g_{m}}_{ac,D_{i,j}^{m}}\leq\frac{1}{\C},\quad\quad  \nn{g^{-1}_{m}}_{ac,g_{m}(D_{i,j}^{m})}\leq\frac{1}{\C}.
   	\end{equation}

    For convenience, we will use the following abbreviations:
    \begin{gather*}
     \hat{B}^{m,\pm}_{i,j}:=g_{m}(\hat{A}^{m,\pm}_{i,j}),\quad\quad
    	C^{m,-}_{i,j}:=\hat{A}^{m,-}_{i,j}, \quad\quad
       C^{m,+}_{i,j}:=\hat{B}^{m,+}_{i,j}\\
      h_{i,j}^{m,-}:=(g_{m})_{\restriction \hat{A}^{m,-}_{i,j}},\quad\quad h_{i,j}^{m,+}:=(g_{m})^{-1}_{\restriction\hat{B}^{m,+}_{i,j}}
´    \end{gather*} 
   
    For $n\geq 0$ and $(i,j)\in\mathcal{P}$ we can find by \cref{f: approximation} and (\ref{e: good approx}) a partition  $\{C_{i,j}^{m,\pm}(1),\dots C_{i,j}^{m,\pm}(\K)\}\subseteq\cl$ of $C^{m,\pm}_{i,j}$ such that:  
    \begin{equation}
    	\label{e: step function approximation}\nn{\mathcal{D}(h^{m,\pm}_{i,j})-\sum_{k=1}^{\K}\frac{k}{\K}\cf_{C^{m,\pm}_{i,j}(k)}}_{1,C^{m,\pm}_{i,j}}\leq\frac{1}{K}.
    \end{equation}
    Up to passing to a further subsequence we may assume that for all $1\leq k\leq \K$ we have
    \begin{equation}
    	\label{e: measures} |\mu(\A_{i,j}^{m,\pm}(k))-\mu(\A_{i,j}^{n,\pm}(k))|,\,\,|\mu(\B_{i,j}^{m,\pm}(k))-\mu(\B_{i,j}^{n,\pm}(k))|\leq\frac{1}{3K^{2}}.
    \end{equation}
    By Observation \ref{o: homoteties} we can find $v,w\in \mh$
    such that for all $(i,j)\in\mathcal{P}$ and $1\leq k\leq \K$ we have 
   
    \begin{enumerate}[label=(\roman{*})]
    	\item \label{idv}$v(\hat{A}^{0,\pm}_{i,j}(k))=\hat{A}^{1,\pm}_{i,j}(k)$,  
    	\item \label{jv}$\mathcal{D}(v)$ has constant value $\frac{\mu(\AA)}{\mu(\aa)}$ on $\aa$,
    	\item \label{idw}$w(\hat{B}^{1,\pm}_{i,j}(k))=\hat{B}^{0,\pm}_{i,j}(k)$,
    	\item \label{jw}$\mathcal{D}(w)$ has constant value $\frac{\mu(\bb)}{\mu(\BB)}$ on $\bb$.
    \end{enumerate} 
    
    Define $A^{+}, A^{-}$ by $A^{\pm}:=\bigcup_{(i,j)\in\{1,\dots r\}^{2}}\bigcup_{k=1}^{\K}\A^{0,\pm}_{i,j}(k)$ and similarly for $B^{\pm}$. From (\ref{e: measures}) and \ref{jv} we get:
      \begin{equation}
      	\label{e: more measures} \nn{v}_{ac,\A^{0,\pm}_{i,j}}=\sum_{k=1}^{\K}|\mu(\A^{0,\pm}_{i,j}(k))-\mu(\A^{1,\pm}_{i,j}(k))|\leq\frac{1}{K},\quad\quad\nn{v}_{ac}\leq\frac{2r^{2}}{K}<\frac{1}{2N}.
      \end{equation}
    Clearly, similar inequalities hold for $w$ and $\B^{l,\pm}_{i,j}$ in place of $v$ and $\A^{l,\pm}_{i,j}$

    Let $\tilde{g}_{1}:=wg_{1}v$ and note that $\tilde{g}_{1}(\aa)=g_{0}(\aa)=\bb$ for all $1\leq i,j\leq r$ and $1\leq k\leq \K$. 
    It is easy to see that there are elements $u,\tilde{u}\in \mh$ with
    $$u_{\restriction A^{+}}=\Id_{A^{+}},\quad\quad\tilde{u}_{\restriction B^{-}}=\Id_{B^{-}},\quad\quad g_{0}=\tilde{u}\tilde{g}_{1}u,$$ and with $u$ preserving each $\aa$ and $\tilde{u}$ preserving each $\bb$.
    
    \begin{lemma}
    	\label{l: final rectification} The inequalities $\nn{u}_{ac},\nn{\tilde{u}}_{ac}<\frac{1}{2N}$ hold.    	
    \end{lemma}
    \begin{subproof}
     We get the following estimate: 
    	\begin{align*}
    		&\nn{\mathcal{D}(\tilde{g}_{1})-\sum_{k=1}^{\K}\frac{k}{\K}\cf_{\hat{A}^{0,-}_{i,j}(k)}}_{1,\hat{A}^{0,-}_{i,j}}=
    		\nn{\mathcal{D}(wg_{1})-\sum_{k=1}^{\K}\frac{k\mu(\hat{A}^{0,-}_{i,j}(k))}{\K\mu(\hat{A}^{1,-}_{i,j}(k))}\cf_{\hat{A}^{1,-}_{i,j}(k)}}_{1,\hat{A}^{1,-}_{i,j}}\leq\\
    		&\nn{\mathcal{D}(g_{1})-\sum_{k=1}^{\K}\frac{k}{\K}\cf_{\hat{A}^{1,-}_{i,j}(k)}}_{1,\hat{A}^{1,-}_{i,j}}+\sum_{k=1}^{\K}\frac{k}{\K}|\mu(\A^{0,-}_{i,j}(k))-\mu(\A^{1,-}_{i,j}(k))|+\nn{\mathcal{D}(wg_{1})-\mathcal{D}(g_{1})}_{1,\hat{A}^{1,-}_{i,j}}\leq\\
    		&\leq\frac{1}{K}+\frac{1}{K}+\nn{w}_{ac,\hat{B}^{1}_{i,j}}\leq\frac{3}{K}\leq\frac{1}{2N}, 
    	\end{align*}
    	where in the equation we use \cref{o: change of variable} as in the proof of \ref{o: change of variable} as well as the identities \ref{idv} and \ref{jv} 
    	and in the last inequality we have used (\ref{e: step function approximation}) to bound the first term, Observation \ref{o: right invariance} once again and \eqref{e: measures} to bound the second term.
    	Adding over all $1\leq i\leq j\leq r$
    	and combining the result with the estimate (\ref{e: step function approximation}) for $m=0$ yields 
    	$\nn{\mathcal{D}(\tilde{g}_{1})-\mathcal{D}(g_{0})}_{1,A^{-}}\leq\frac{4(|\mathcal{P}^{+}|+|\mathcal{P}^{-}|)}{K}=\frac{1}{2N}$. 
    	Using Observation \ref{o: right invariance} again and then the inequality above we obtain:    
    	$$
    	 \nn{\tilde{u}}_{ac}=\nn{\tilde{u}}_{ac,B^{-}}=\nn{\mathcal{D}(\tilde{u}\tilde{g}_{1})-\mathcal{D}(\tilde{g}_{1})}_{1,B^{-}}=\nn{\mathcal{D}(g_{0})-\mathcal{D}(\tilde{g}_{1})}_{1,B^{-}}\leq\frac{1}{2N}
    	$$
    	By estimating $\nn{\mathcal{D}(\tilde{g}_{1}^{-1})-\mathcal{D}(g_{0}^{-1})}_{1,B^{+}}$ symmetrically we also get $\nn{u}_{ac}\leq\frac{1}{2N}$.      	
    \end{subproof}  
   Now $g_{0}=\tilde{u}wg_{1}vu$ with $\tilde{u}w,vu\in \mh$. On the other hand, $\nn{\tilde{u}w}_{ac}\leq\nn{\tilde{u}}_{ac}+\nn{w}_{ac}\leq\frac{1}{N}$ by \ref{o: right invariance} and similarly $\nn{vu}\leq\frac{1}{N}$, so $\tilde{u}w,vu\in\mathcal{W}$, contradicting the choice of $(g_{n})_{n\geq 0}$. This concludes the proof. 
   \end{proof}

   \bibliographystyle{plain}
  \bibliography{bibliography}

\begin{thebibliography}{10}

\bibitem{anderson1958algebraic}
Richard~D Anderson.
\newblock The algebraic simplicity of certain groups of homeomorphisms.
\newblock {\em American Journal of Mathematics}, 80(4):955--963, 1958.

\bibitem{bartle2001modern}
Robert~Gardner Bartle.
\newblock {\em A modern theory of integration}, volume~32.
\newblock American Mathematical Society Providence, 2001.

\bibitem{choksi1958inverse}
JR~Choksi.
\newblock Inverse limits of measure spaces.
\newblock {\em Proceedings of the London Mathematical Society}, 3(3):321--342,
  1958.

\bibitem{CV1977}
Charles~O. Christenson and William~L. Voxman.
\newblock {\em Aspects of topology}.
\newblock Pure and Applied Mathematics, Vol. 39. Marcel Dekker, Inc., New
  York-Basel, 1977.

\bibitem{cohen2017polishability}
Michael~P Cohen.
\newblock Polishability of some groups of interval and circle diffeomorphisms.
\newblock {\em arXiv preprint arXiv:1709.04523}, 2017.

\bibitem{edwards1971deformations}
Robert~D Edwards and Robion~C Kirby.
\newblock Deformations of spaces of imbeddings.
\newblock {\em Annals of Mathematics}, 93(1):63--88, 1971.

\bibitem{fathi1980structure}
Albert Fathi.
\newblock Structure of the group of homeomorphisms preserving a good measure on
  a compact manifold.
\newblock In {\em Annales scientifiques de l'{\'E}cole Normale Sup{\'e}rieure},
  volume~13, pages 45--93, 1980.

\bibitem{gartside2003autohomeomorphism}
Paul Gartside and Aneirin Glyn.
\newblock Autohomeomorphism groups.
\newblock {\em Topology and its Applications}, 129(2):103--110, 2003.

\bibitem{herndon2018absolute}
Jake Herndon.
\newblock Absolute continuity and large-scale geometry of polish groups.
\newblock {\em arXiv preprint arXiv:1802.10239}, 2018.

\bibitem{ihli2022generation}
Dakota~Thor Ihli.
\newblock Generation and genericity of the group of absolutely continuous
  homeomorphisms of the interval.
\newblock {\em arXiv preprint arXiv:2203.02854}, 2022.

\bibitem{kechris2007turbulence}
Alexander~S Kechris and Christian Rosendal.
\newblock Turbulence, amalgamation, and generic automorphisms of homogeneous
  structures.
\newblock {\em Proceedings of the London Mathematical Society}, 94(2):302--350,
  2007.

\bibitem{kirby2010stable}
Robion Kirby.
\newblock Stable homeomorphisms and the annulus conjecture.
\newblock In {\em Topological Library: Part 2: Characteristic Classes and
  Smooth Structures on Manifolds}, pages 253--261. World Scientific, 2010.

\bibitem{klenke2013probability}
Achim Klenke.
\newblock {\em Probability theory: a comprehensive course}.
\newblock Springer Science \& Business Media, 2013.

\bibitem{mann2016automatic}
Kathryn Mann.
\newblock Automatic continuity for homeomorphism groups and applications.
\newblock {\em Geometry \& Topology}, 20(5):3033--3056, 2016.

\bibitem{oxtoby1941measure}
John~C Oxtoby and Stanislas~M Ulam.
\newblock Measure-preserving homeomorphisms and metrical transitivity.
\newblock {\em Annals of Mathematics}, pages 874--920, 1941.

\bibitem{quinn1982ends}
Frank Quinn.
\newblock Ends of maps. iii. dimensions 4 and 5.
\newblock {\em Journal of Differential Geometry}, 17(3):503--521, 1982.

\bibitem{rado2012continuous}
Tibor Rad{\'o} and Paul~V Reichelderfer.
\newblock {\em Continuous transformations in analysis: with an introduction to
  algebraic topology}, volume~75.
\newblock Springer Science \& Business Media, 2012.

\bibitem{rosendal2013global}
Christian Rosendal.
\newblock Global and local boundedness of polish groups.
\newblock {\em Indiana University Mathematics Journal}, pages 1621--1678, 2013.

\bibitem{rosendal2021coarse}
Christian Rosendal.
\newblock {\em Coarse geometry of topological groups}, volume 223.
\newblock Cambridge University Press, 2021.

\bibitem{rosendal2006automatic}
Christian Rosendal and Slawomir Solecki.
\newblock Automatic continuity of homomorphisms and fixed points on metric
  compacta.
\newblock {\em Israel Journal of Mathematics}, 162:349–371, 2007.

\bibitem{solecki11999polish}
S{\l}awomir Solecki.
\newblock Polish group topologies.
\newblock In S~Barry Cooper and John~K Truss, editors, {\em Sets and proofs},
  volume 258. Cambridge University Press, 1999.

\end{thebibliography}

\end{document}